%% file: bigpw.tex
\documentclass[12pt]{article}
\usepackage{amsmath}
\usepackage{amsthm}
\usepackage{amssymb}
\usepackage{graphicx}
\usepackage{epsfig}
\usepackage{pgf,tikz}
\usetikzlibrary{arrows}
\usepackage[english]{babel}
\usepackage[utf8]{inputenc}
\usepackage{tikz}
\usepackage{verbatim}
\usepackage[colorlinks=false]{hyperref}
\usepackage[margin=1.1in]{geometry}
\usepackage{soul}
\usepackage[normalem]{ulem}
\usepackage{xcolor}

\def\he{homeomorphic embedding}
\def\emb{\hookrightarrow}
\def\td{tree-de\-com\-position}
\def\pd{path-de\-com\-position}
\def\sR{{\cal R}}

\def\dfn#1{{\em #1}}
\def\claimm#1#2\par{\smallbreak\noindent\rlap{\rm(#1)}\ignorespaces
\hangindent=30pt\hskip30pt{\ignorespaces\sl#2}\smallskip}

\title{Minors of 2-connected graphs \ca{of} high path-width}
\author{Thanh N. Dang  \hspace{20 mm} Robin Thomas\\ \\
Georgia Institute of Technology}

\begin{document}
\centerline{\Large \bf MINORS OF TWO-CONNECTED GRAPHS}
\medskip
\centerline{{\Large\bf OF LARGE PATH-WIDTH}%
\footnote{Partially supported by NSF under Grants No.~DMS-1202640 and~DMS-1700157. 11 October 2015, revised 14 April 2018.
}}

\bigskip
\bigskip

\centerline{{\bf Thanh N. Dang}%
}
\smallskip
\centerline{and}
\smallskip
\centerline{{\bf Robin Thomas}}
\bigskip
\centerline{School of Mathematics}
\centerline{Georgia Institute of Technology}
\centerline{Atlanta, Georgia  30332-0160, USA}
\bigskip

\newtheorem{theorem}{Theorem}[section]
\newtheorem{lemma}[theorem]{Lemma}
\newtheorem{conjecture}{Conjecture}
\newtheorem{problem}{Open problem}
\newtheorem{proposition}[theorem]{Proposition}
\newtheorem{corollary}[theorem]{Corollary}
\newtheorem{claim}{Claim}[theorem]
\newenvironment{definition}[1][Definition]{\begin{trivlist}
\item[\hskip \labelsep {\bfseries #1}]}{\end{trivlist}}
\newenvironment{example}[1][Example]{\begin{trivlist}
\item[\hskip \labelsep {\bfseries #1}]}{\end{trivlist}}
\newenvironment{remark}[1][Remark]{\begin{trivlist}
\item[\hskip \labelsep {\bfseries #1}]}{\end{trivlist}}
\def\ca#1{{\color{black}#1}}
\def\cb#1{{\color{red}#1}}
\def\cc#1{{\color{blue}#1}}
\def\cd#1{{\color{violet}#1}}
 \def\cb#1{#1}
\def\cc#1{#1}
\def\cd#1{#1}
\newcommand\redsout{\bgroup\markoverwith{\textcolor{red}{\rule[0.5ex]{2pt}{1pt}}}\ULon}
 \def\redsout#1{}

\begin{abstract}
\noindent
Let $P$ be a graph with a vertex $v$ such that $P\backslash v$ is a forest, and let $Q$ be an outerplanar graph.
We prove that there exists a number $p=p(P,Q)$ such that every
$2$-connected graph of path-width at least $p$ has a minor isomorphic to $P$ or~$Q$. This result answers
a question of Seymour and implies a conjecture of Marshall and Wood.
The proof is based on a new  property of \td s.
\end{abstract}

\section{Introduction}
All {\em graphs} in this paper are finite and simple; that is, they  have no loops or parallel edges.
 {\em Paths} and {\em cycles} have no ``repeated" vertices or edges.
A graph $H$ is a {\em minor} of a graph $G$ if we can obtain $H$ by contracting edges of a subgraph of $G$.
An {\em$H$ minor} is a minor isomorphic to $H$.
A tree-decomposition of a graph $G$ is a pair $(T,X)$, where $T$ is a tree and $X$ is a family
\ca{$(X_t : t \in V(T))$} such that:

\begin{itemize}
\item[(W1)]
 $ \bigcup_{t\in\ V((T)} X_t = V(G)$, and for every edge of $G$ with ends $u$ and $v$
there exists $t \in V(T)$ such that $u,v \in X_t$, and

\item[(W2)] if $\ca{t_1}, t_2, t_3 \in V(T)$ and $t_2$ lies on the path in $T$ between $t_1$ and $t_3$,
then $X_{t_1} \cap X_{t_3} \subseteq X_{t_2}$.
\end{itemize}
The {\em width} of a tree-decomposition $(T,X)$ is $\max\{ |X_t| - 1:t\in V(T)\}$.
The {\em tree-width} of a graph $G$ is the smallest width among all tree-decompositions of $G$.
A {\em path-decomposition} of $G$ is a tree-decomposition $(T,X)$ of $G$, where $T$ is a path.
We will often denote a \pd\ as $(X_1,X_2,\ldots,X_n)$, rather than having the constituent sets
indexed by the vertices of a path.
The {\em path-width} of $G$ is the smallest width among all path-decompositions of $G$.
\ca{Robertson and Seymour \cite{gm5}} proved the following:
\begin{theorem}
\label{thm1}
For every  planar graph $H$ there exists an integer $n = n(H)$ such that  every graph of tree-width at least $n$ has an $H$  minor.
\end{theorem}
Robertson and Seymour \cite{robertson83} \ca{also} proved an analogous result for path-width:
\begin{theorem}
\label{thm2}
For every forest $F$, there exists an integer $p=p(F)$ such that every graph of path-width at least $p$ has an $F$  minor.
\end{theorem}
\noindent
Bienstock,  Robertson, Seymour and the second author~\cite{bienstock912} gave a simpler proof of Theorem~\ref{thm2} and
improved the value of $p$ to $|V(F)| - 1$, which  is  best possible,  because  $K_k$  has path-width $k-1$ and does not have any
forest minor on $k+1$ vertices. A yet simpler proof of Theorem~\ref{thm2} was found by Diestel~\cite{DieGM1}.

While Geelen, Gerards and Whittle~\cite{geelen} generalized Theorem~\ref{thm1} to representable matroids,
it is not {\em a priori} clear what a version  of Theorem~\ref{thm2} for matroids should be,  because excluding
a forest in matroid setting is equivalent to imposing a bound on the number of elements and has  no relevance to path-width. To overcome this, Seymour
\cite[Open Problem~2.1]{dean} asked if there was a generalization of Theorem~\ref{thm2} for $2$-connected graphs with forests replaced by the two
families of graphs  mentioned in the abstract.
Our main result answers Seymour's question in the affirmative:
\begin{theorem}
\label{thm3}
Let $P$ be a graph with a vertex $v$ such that $P\backslash v$ is a forest,
and let $Q$ be an outerplanar graph. Then there exists a number $p=p(P,Q)$
such that every $2$-connected graph of path-width at least $p$
has a $P$ or $Q$  minor.
\end{theorem}

Theorem~\ref{thm3} is a generalization of Theorem~\ref{thm2}.
To deduce Theorem~\ref{thm2} from Theorem~\ref{thm3},
 given a graph $G$, we may assume that $G$ is connected, because the path-width of a graph is equal
to the maximum path-width of its components. We add one vertex and make it adjacent to every vertex of $G$.
Then the new graph is $2$-connected, and by Theorem~\ref{thm3},
it has a $P$ or $Q$  minor. By choosing suitable $P$ and $Q$,
we can get an $F$ minor in $G$.

Our strategy to prove Theorem~\ref{thm3} is as follows. Let $G$ be a $2$-connected graph of large path-width.
We may assume that the tree-width of $G$ is bounded, for otherwise $G$ has a minor isomorphic to both $P$ and $Q$
by Theorem~\ref{thm1}. So let $(T,X)$ be a tree-decomposition of $G$ of bounded width. Since the path-width of $G$
is large, it follows by a simple argument (Lemma~\ref{lemma4.3} below) that the path-width of $T$ is large, and hence
it has a subgraph $T'$ isomorphic to a subdivision of a  large binary tree by Theorem~\ref{thm2}.
It now seems plausible that we could use $T'$ and  properties (W3) and (W4) of tree-decompositions, introduced below, 
which we can assume by~\cite{OpoOxlTho,thomas90}, to show the desired conclusion.
But there is a catch: for instance, a long cycle has a \td\ $(T,X)$ satisfying (W3) and (W4) (and, in fact, the minimality
condition used in their proof, as well as that of Bellenbaum and Diestel~\cite{BelDie}) such that $T$ has a subgraph
isomorphic to a large binary tree.
And yet it feels that this is the ``wrong" \td\ and that the ``right" \td\ is one where $T$ is a path.
The main result of the first part of this paper, Theorem~\ref{w7} below, deals with
converting these ``branching" \td s into ``non-branching" ones without increasing their width.

Marshall and Wood \cite{marshall}
define $g(H)$ as the minimum number for which there exists a positive integer
$p(H)$ such that every $g(H)$-connected graph with no $H$ minor has path-width
at most $p(H)$. Then Theorem~\ref{thm2} implies that $g(H) = 0$ if and only if $H$ is a forest.
There is no graph $H$ with $g(H) = 1$, because path-width of a graph $G$
is the maximum of the \ca{path-widths} of its connected components.
Let $A$ be the graph that consists of a cycle $a_1a_2a_3a_4a_5a_6a_1$
and extra edges $a_1a_3, a_3a_5, a_5a_1$. Let $C_{3,2}$ be the graph consisting of two disjoint triangles.
In Section~\ref{sec2} we prove a conjecture of
Marshall and Wood \cite{marshall}:
\begin{theorem}
\label{con1}
A graph $H$ has no $K_{4}, K_{2,3}$, $C_{3,2}$ or $A$  minor if  and only if $g(H) \leq 2 $.
\end{theorem}
In Section~\ref{sec:linked} we describe a special tree-decomposition, whose existence we establish in Section~\ref{sec:thm}.
Section~\ref{sec:quasi}  introduces a quasi-order on  trees, our main tool in the proof of Theorem~\ref{w7}. 
In Section~\ref{sec:cascades} we introduce ``cascades", our main tool in the proof 
of Theorem~\ref{thm3}, and prove that in any \td\
with no duplicate bags of  bounded width of a graph of big path-width there is an ``injective" cascade
of large height. In Section~\ref{sec:ordered} we prove that every $2$-connected graph of big path-width
and bounded tree-width admits a \td\ of bounded width and a cascade with linkages that are minimal.
In Section~\ref{sec:taming} we analyze those minimal linkages and prove that there are essentially only two
types of minimal linkage. This is where we use the properties of \td s from Section~\ref{sec:linked}.
Finally,  in Section~\ref{sec4} we convert the two types of linkage into the two families of graphs
from
Theorem~\ref{thm3}.

\section{Proof of Theorem~\ref{con1}}
\label{sec2}
In this section we prove that Theorem~\ref{con1} is implied by Theorem~\ref{thm3}.

\begin{definition}
Let $h\ge0$ be an integer.
By a {\em binary tree of height $h$}  we mean a  tree with a unique vertex $r$ of degree two
and all other vertices of degree one or three such that every vertex of degree one is at distance
exactly $h$  from $r$.
Such a tree is unique up to isomorphism and so we will speak of the binary tree of height $h$.
We denote the binary tree of height $h$ by $CT_h$ and we call $r$ the {\em root} of $CT_h$.
Each vertex in $CT_h$ with distance $k$ from $r$ has {\em height} $k$.
We call  the vertices at distance $h$ from $r$ the {\em leaves of $CT_h$.}
If $t$ belongs to the unique path in $CT_h$ from $r$
to a vertex $t'\in V(T_h)$, then we say
that $t'$  is a {\em  descendant} of $t$ and that $t$
is an {\em  ancestor} of $t'$.
If, moreover, $t$ and $t'$ are adjacent, then we say that $t$ is the {\em parent}
of $t'$ and that $t'$ is a {\em child} of $t$.

Let $\ca{\mathcal{P}_k}$ be the graph consisting of $CT_k$ and a separate vertex that is adjacent to every leaf of $CT_k$.
\end{definition}

\begin{lemma}
\label{lem:nearforest}
If a graph $H$ has no $K_4, C_{3,2}$, or $A$  minor,
then   $H$ has a vertex $v$ such that $H\backslash v$ is a forest.
\end{lemma}

\noindent\textit{Proof.}
We proceed by induction on $|VH)|$.
The lemma clearly holds when $|V(H)|=0$, and so we may assume that $H$ has at least one vertex and
that the lemma holds for graphs on fewer than $|V(H)|$ vertices.
If $H$ has a vertex of degree at most one, then the lemma follows by induction by deleting such  vertex.
We  may therefore assume that $H$ has minimum degree at least two.

If  $H$  has a cutvertex, say $v$,  then $v$ is as desired,  for if $C$ is a cycle in $H\backslash v$,
then $H\backslash V(C)$ also contains a cycle (because $H$ has minimum degree at least two), and hence
$H$ has a $C_{3,2}$ minor, a  contradiction.  We may therefore assume that $H$  is $2$-connected.

We may assume that $H$ is not a cycle, and hence it has
an ear-decomposition $H=H_0\cup H_1\cup\cdots\cup H_k$, where $k\ge1$, $H_0$ is a cycle and
for $i=1,2,\ldots,k$ the graph $H_i$ is a path with ends $u_i,v_i\in V(H_0\cup H_1\cup\cdots\cup H_{i-1})$
and otherwise disjoint from $H_0\cup H_1\cup\cdots\cup H_{i-1}$.
If $u_1\in\{u_i,v_i\}$ for all $i\in\{2,3,\ldots,k\}$, then $u_1$  satisfies the conclusion of the lemma, and similarly for $v_1$.
We may therefore assume that there exist $i,j\in \{2,3,\ldots,k\}$ such that $u_1\not\in\{u_i,v_i\}$ and
$v_1\not\in\{u_j,v_j\}$.
It  follows that $H$  has a $K_4, C_{3,2}$, or $A$  minor, a contradiction.~\qed

\begin{lemma}
\label{lemma2.1}
If a graph  $H$ has a vertex $v$ such that $H\backslash v$ is a forest.
then there exists an integer $k$ such that  $H$ is isomorphic to a minor of $\ca{\mathcal{P}_k}$.
\end{lemma}

\noindent\textit{Proof.}
Let $v$ be such that $T:=H\backslash v$ is a forest. We may assume, by replacing $H$ by a graph with an
$H$ minor, that $T$ is isomorphic to $CT_t$ for some $t$, and that $v$ is adjacent to every vertex of $T$.
It follows that $H$ is isomorphic to a minor of $\mathcal{P}_{2t}$, as desired.~\qed

\begin{definition}
Let $\ca{\mathcal{Q}_1}$ be $K_3$. An arbitrary edge of $\mathcal{Q}_1$ will be designated as {\em base edge}.
For $i\ge2$ the graph $\ca{\mathcal{Q}_{i}}$ is constructed as follows:
Now assume that $\mathcal{Q}_{i-1}$ has already been defined, and let $Q_1$ and $Q_2$ be two disjoint copies of	$\mathcal{Q}_{i-1}$
with base edges $u_1v_1$ and $u_2v_2$, respectively.
Let $T$ be a copy of $K_3$ with vertex-set $\{w_1,w_2,w\}$ disjoint from $Q_1$ and $Q_2$.
The graph $\mathcal{Q}_{i}$ is obtained from $Q_1\cup Q_2\cup T$ by identifying $u_1$ with $w_1$,
 $u_2$ with $w_2$, and $v_1$ and $v_2$ with $w$.
The edge $w_1w_2$ will be the {\em base edge} of  $\mathcal{Q}_{i}$.

A graph is {\em outerplanar} if it has a drawing in the plane (without crossings) such that every vertex is incident
with the unbounded face.
A graph is a {\em near-triangulation} if it is drawn in the plane in such a way that every face except possibly
the unbounded one is bounded by a triangle.

Let $H$ and $G$ be graphs.  If $G$ has an $H$ minor, then to every vertex $u$ of $H$ there corresponds a connected
subgraph of $G$, called the {\em node of $u$}.
\end{definition}

\begin{lemma}
\label{lem:neartriang}
Let $H$ be a $2$-connected outerplanar near-triangulation with $k$ triangles. Then
 $H$ is isomorphic to a minor of $\ca{\mathcal{Q}_k}$.
Furthermore, the minor inclusion can be chosen in such a way that for every edge $a_1a_2\in E(H)$ incident with the unbounded face
and  for every $i\in\{1,2\}$, the vertex $w_i$ belongs to the node of $a_i$,
where $w_1w_2$ is the base edge of~$\mathcal{Q}_k$.
\end{lemma}

\begin{proof}
We proceed by induction on $k$.
The lemma clearly holds when $k=1$, and so we may assume that $H$ has at least two triangles and
that the lemma holds for graphs with fewer than $k$ triangles.
The edge $a_1a_2$ belongs to a unique triangle, say $a_1a_2c$. The triangle $a_1a_2c$ divides $H$ into two near-triangulations
$H_1$ and  $H_2$, where the edge $a_ic$ is incident with the unbounded face of $H_i$.
Let $Q_1,Q_2,u_1,v_1,u_2,v_2,w_1,w_2$ be as in the definition of~$\mathcal{Q}_k$.
By the induction hypothesis the graph $H_i$ is isomorphic to a minor of $Q_i$ in such a way that the vertex $u_i$ belongs to
the node of $a_i$ and the vertex $v_i$ belongs to the node of $c$.
It follows that $H$ is isomorphic to $\mathcal{Q}_k$ in such a way that $w_i$ belongs to the node of~$a_i$.
\end{proof}

\begin{lemma}
\label{lemma2.2}
Let $H$ be a graph that has no $K_4$ or $K_{2,3}$  minor.
Then there exists an integer $k$ such that
 $H$ is isomorphic to a minor of $\ca{\mathcal{Q}_k}$.
\end{lemma}

\begin{proof}
It is well-known~\cite[Exercise~23]{DieGT5} that the hypotheses of the lemma imply that $H$
is outerplanar.
 We may assume, by replacing $H$ by a graph with an
$H$ minor, that $H$ is a $2$-connected outerplanar near-triangulation.
The lemma now follows from Lemma~\ref{lem:neartriang}.
\end{proof}

\begin{corollary}
\label{lemma2.3}
Let $H$ be a graph that has no $K_4$, $K_{2,3}$, $C_{3,2}$,
or $A$  minor. Then there exists an integer $k$ such that
 $H$ is isomorphic to a minor of $\ca{\mathcal{P}_k}$ and $H$ is isomorphic to a minor of $\ca{\mathcal{Q}_k}$.
\end{corollary}

\begin{proof}
This follows from Lemmas~\ref{lem:nearforest}, \ref{lemma2.1} and~\ref{lemma2.2}.
\end{proof}

\begin{proof}[Proof of Theorem~\ref{con1}, assuming Theorem~\ref{thm3}]
To prove the``if" part
notice that $\ca{\mathcal{P}_k}$ and $\ca{\mathcal{Q}_k}$ are $2$-connected and have
large path-width when $k$ is large,
 because $\mathcal{Q}_k$ has a $CT_{k-1}$ minor.
There is no vertex $v$ in $A$ such that $A\backslash v$ is acyclic.
So, $A$ and $C_{3,2}$ are not minors of $\ca{\mathcal{P}_k}$ for any $k$. The graph
$\ca{\mathcal{Q}_k}$ is outerplanar, so $K_4$ and $K_{2,3}$ are not minors of
$\ca{\mathcal{Q}_k}$ for any positive integer $k$. This means $g(H) \geq 3$ for
$H \in \{K_4, K_{2,3}, C_{3,2}, A\}$. This proves the ``if" part.

To prove the ``only if" part,
if $H$ has no $K_{4},
K_{2,3}$, $C_{3,2}$ or $A$ minor, then by Corollary~\ref{lemma2.3}
$H$ is a minor of both $\ca{\mathcal{P}_k}$ and $\ca{\mathcal{Q}_k}$ for some $k$.
Then $g(H) \leq 2$ by Theorem~\ref{thm3}.
\end{proof}

\input w7insert.tex

\section{Cascades}
\label{sec:cascades}
In this section we introduce ``cascades", our main tool. The main result of this section,
Lemma~\ref{lemma4.8}, states that in any \td\
with no duplicate bags of  bounded width of a graph of big path-width there is an ``injective" cascade
of large height

\begin{lemma}
\label{lemma4.3}
Let $p, w$ be two positive integers and let $G$ be a graph of tree-width strictly  less than
$w$ and path-width at least $p$. Then for every tree-decomposition $(T,X)$ of $G$ of
width strictly less than $w$, the path-width of $T$ is at least $\lfloor p/w \rfloor$ .
\end{lemma}

\begin{proof}
We will prove the contrapositive.
Assume there exists a tree-decomposition $(T,X)$ of $G$ of width $< w$ such that
the path-width of $T$ is less than $\lfloor p/w \rfloor$.
Because the path-width of $T$ is less than $\lfloor p/w \rfloor$, there exists a path-decomposition
$(Y_1, Y_2,...,Y_s)$ of $T$ with $|Y_i| \leq \lfloor p/w \rfloor$ for all $i$.
We will construct a path-decomposition $ (Z_1, Z_2, ..., Z_s)$ for
$G$ of width less than $p$.
Set $Z_i = \bigcup_{y \in Y_i} X_{y}$ for every $i \in \{1,2,...,s\}$.
For every vertex $v \in V(G)$, $v$ belongs
to at least one set $X_t$ for some $t \in V(T)$. The vertex $t$ of the tree $T$ must be in $Y_l$
for some $l \in \{1,2,...,s\}$,
so $v \in X_t \subseteq Z_l$. Therefore, $\bigcup Z_i = V(G)$.
Similarly, for every edge $uv \in E(G)$, there exists $t\in V(T)$ such that
$u,v \in X_t$. Therefore, $u, v\in Z_l$ for some $l \in \{1,2,...,s\}$.

Now, if a vertex $v \in V(G)$ belongs to both $Z_a$ and $Z_b$ for some
$a, b\in \{1,2,...,s\}, a<b$, we will show that
$v \in Z_c$ for all $c$ such that $a < c < b$. Let $c$ be an arbitrary integer satisfying $a<c<b$.
 \ca{The fact that} $v \in Z_a$ implies $v \in X_{y_1}$
for some $y_1 \in Y_a$. Similarly, $v \in X_{y_2}$ for some $y_2 \in Y_b$.
Let $H$ be the set of vertices of $T$ on the path from $y_1$ to $y_2$.
Since $y_1 \in Y_a$ and $y_2 \in Y_b$, $H\cap Y_{a} \neq \emptyset \neq H\cap Y_{b}$.
Hence, by Lemma~\ref{cutbag} with $H = T$ and $(T,Y)$ \ca{the path-decomposition} $(Y_1, Y_2,...,Y_s)$,
we have $H \cap Y_{c} \neq \emptyset$.
Let $t\in H \cap Y_{c}$, then $v \in X_t \subseteq Z_c$.
So $(Z_1, Z_2, ..., Z_s)$ is a path-decomposition of $G$.
Since \ca{the} width of $(T,X)$ is less than $w$,
we have $|X_{y}| \leq w$ for every $y \in Y_i$, where $i\in \{1,2,...,s\}$.
Therefore, $|Z_i| \leq w.\lfloor p/w \rfloor \leq p $ for every $i\in\{1,2,...,s\}$.
Therefore, the width of $(Z_1, Z_2, ..., Z_s)$ is less than $p$,
so the path-width of $G$ \ca{is} less than $p$, \ca{as desired}.
\end{proof}

Let $T,T'$ be trees. A {\em \he\ of $T$  into $T'$}  is a mapping $\eta:V(T)\to V(T')$
such that
\begin{itemize}
\item $\eta$ is an injection, and
\item if $tt_1,tt_2$ are edges of $T$  with a common end, and $P_i$ is the unique path in $T'$
 with ends $\eta(t)$ and $\eta(t_i)$, then $P_1$ and $P_2$ are edge-disjoint.
\end{itemize}
We will write $\eta:T\emb T'$ to denote that $\eta$ is a \he\ of $T$  into $T'$.
Since $CT_a$ has maximum degree at most three, the following lemma follows
from~\cite[Lemma~6]{marshall}.
\begin{lemma}
\label{lemma4.4}
Let $T$ be a forest of path-width at least $a \geq 1$.
Then there exists a \he\ $CT_{a-1}\emb T$.
\end{lemma}

For every integer $h\ge1$ we will need a specific type of tree,
which we will denote by $T_h$.
The  tree $T_h$  is obtained from $CT_h$ by subdividing every edge
not incident with a vertex of degree one exactly once, and adding a new vertex $r'$ of degree one adjacent
to the root $r$ of $CT_h$.
The vertices of $T_h$  of degree three will be called {\em  major},
and all the other vertices
will be called {\em  minor}.
We say that $r$ is the  {\em major root} of $T_h$  and that $r'$ is the
{\em minor root} of $T_h$.
Each major vertex at distance $2k$ from $r$ has {\em height} $k$,
and each minor vertex at distance $2k$ from $r'$ has {\em height} $k$.

If $t$ belongs to the unique path in $T_h$ from $r'$
to a vertex $t'\in V(T_h)$, then we say
that $t'$  is a {\em  descendant} of $t$ and that $t$
is an {\em  ancestor} of $t'$.
If, moreover, $t$ and $t'$ are adjacent, then we say that $t$ is the {\em parent}
of $t'$ and that $t'$ is a {\em child} of $t$.
Thus every major vertex $t$ has exactly three minor neighbors.
Exactly one of those neighbors is an ancestor of $t$.
The other two neighbors are descendants
of $t$.  We will assume that one of the two  descendant
neighbors is designated as the {\em  left neighbor}
and  the other as the {\em  right neighbor}.
Let $t_0,t_1,t_2$  be the parent, left neighbor and right neighbor of $t$, respectively.
We say that the ordered triple $(t_0,t_1,t_2)$ is the
{\em trinity at $t$}. In case we want to emphasize that
the trinity is at $t$, we use the notation $(t_0 (t) ,t_1(t),t_2(t))$.

Let $\eta:T\emb T'$ . We define
$sp(\eta)$, the {\em span of $\eta$}, to be the set of vertices $t \in V(T') $ that lie on the path
from $\eta(t_1)$ to $\eta(t_2)$ for some
vertices $t_1, t_2 \in V(T)$.

Let $s > 0$  be an integer and  let $(T,X)$  be  a tree-decomposition of a graph $G$.
By a {\em cascade of height $h$ and size $s$ in $(T,X)$} we mean a \he\ $\eta:T_h\emb T$ such that
$|X_{\eta(t)}|=s$ for every minor vertex $t\in V(T_h)$ and $|X_{t}|\geq s$ for every
$t$ in the span of $\eta$.

\begin{lemma}
\label{lemma4.5}
For any positive integer $h$ and nonnegative integers $a,k$, the following holds.
Let $m = (a+2)h + a$.
Let $(T, X)$ be a tree-decomposition of a graph $G$
and let $\phi: CT_{m} \emb T$ be a homeomorphic embedding such that
$|X_t| \geq k$ for all $t \in sp(\phi)$.
If for every $t\in V(CT_m)$ at height  $l\le m-a$ there exist a descendant $t'$ of $t$
at height $l+a$ and a vertex $r\in T[\phi(t),\phi(t')]$
 such that $|X_{r}|=k$,
then there exists a cascade $\eta$ of height $h$ and size $k$ in $(T, X)$.
\end{lemma}

\begin{proof}
By hypothesis there exist a vertex $x_0\in V(CT_m)$ at height  $a$ and
a vertex $u_0\in V(T)$ on the path from the image under $\phi$ of the root
of $CT_m$ to $\phi(x_0)$
such that $|X_{u_0}|=k$.
Let $x$ be a child of $x_0$, and let $x_1$ and $x_2$
be the children of $x$. By hypothesis there exist, for $i=1,2$, a vertex  $y_i\in V(CT_m)$
at height  $2a+2$  that  is a descendant of $x_i$ and a vertex $u_i\in T[\phi(x_i),\phi(y_i)]$
such that  $|X_{u_i}|=k$.
Let $r$ be the major root of $T_1$, and let $(t_0,t_1,t_2)$ be its trinity.
We define $\eta_1:T_1\emb T$ by $\eta_1(t_i)=u_i$  for $i=0,1,2$ and $\eta_1(r)=\phi(x)$.
Then $\eta_1$ is a cascade of  height one and size $k$ in $(T,X)$.
If $h=1$, then $\eta_1$ is as desired, and so we may assume that $h>1$.

Assume now that for some positive integer $l<h$ we have constructed a cascade $\eta_l:T_l\emb T$
of height $l$ and size $k$ in $(T, X)$ such that for every leaf $t_0$ of $T_l$ other than the minor root
there exists a vertex $x_0\in V(CT_m)$ at height $(a+2)l+a$ such that the image under $\eta_l$
of every vertex on the path in $T_l$ from the minor root to $t_0$ belongs to the path in $T$
from the image under $\phi$ of the root of $CT_m$ to $\phi(x_0)$.
Our objective is to extend $\eta_l$ to a cascade $\eta_{l+1}$ of height $l+1$ and size $k$ in $(T, X)$
with the same property.
To that end let $\eta_{l+1}(t)=\eta_l(t)$ for all $t\in V(T_l)$, let $t_0$ be a leaf of $T_l$
other than the minor root and let $x_0$ be as earlier in the paragraph.
Let $x$ be a child of $x_0$, and let $x_1$ and $x_2$
be the children of $x$. By hypothesis there exist, for $i=1,2$, a vertex  $y_i\in V(CT_m)$
at height  $(a+2)(l+1)+a$  that  is a descendant of $x_i$ and a vertex $u_i\in T[\phi(x_i),\phi(y_i)]$
such that  $|X_{u_i}|=k$.
Let $r$ be the  child of $t_0$ in $T_{l+1}$, and let $(t_0,t_1,t_2)$ be its trinity.
We define  $\eta_{l+1}(t_i)=u_i$  for $i=1,2$ and $\eta_{l+1}(r)=\phi(x)$.
This completes the definition of $\eta_{l+1}$.

Now $\eta_h$ is as desired.
%
\end{proof}

\begin{lemma}
\label{lemma4.6}
For any two positive integers $h$ and $w$,
there exists a positive integer $p = p(h, w)$ such that
if $G$ is a graph
of path-width at least $p$,
then in any tree-decomposition
of $G$ of width less than $w$,
there exists a cascade of height $h$.
\end{lemma}

\begin{proof}
%
Let $a_{w+1}= 0$, and for $k=w,w-1,\ldots,0$ let $a_k = (a_{k+1} + 2)h + a_{k+1}$, and let
$p=w(a_0+1)$. We claim that $p$ satisfies the conclusion of the lemma.
To see that let $(T,X)$ be a tree-decomposition of $G$ of width less than $w$.
Let $k\in\{0,1,\ldots,w+1\}$ be the maximum integer such that
there exists a homeomorphic embedding $\phi:CT_{a_{k}}\emb T$
satisfying $|X_{t}|\ge k$ for all $t\in sp(\phi)$.
Such an integer exists, because $k=0$ satisfies those requirements by  Lemmas~\ref{lemma4.3}
and~\ref{lemma4.4}, and it satisfies $k\le w$, because the width of $(T,X)$ is less than $w$.
The maximality of $k$ implies that for the integers $h,k$ and $a_{k+1}$ the hypothesis of Lemma~\ref{lemma4.5}
is satisfied.
Thus the lemma follows from Lemma~\ref{lemma4.5}.
%
\end{proof}

 Let $(T,X)$  be  a tree-decomposition of a graph $G$, and let  $\eta:T_h\emb T$ be a
cascade of height $h$ and size $s$ in $(T,X)$.
We say that $\eta$  is {\em injective} if there exists $I \subseteq V(G)$ such that
$|I| < s$ and $X_{\eta(t)}\cap X_{\eta(t')}= I$
for every two distinct vertices $t,t'\in V(T_h)$. We call this set $I$ the
{\em common intersection set} of $\eta$.

\begin{lemma}
\label{lemma4.7}
 Let $a, b, s,w$ be positive integers and let
$k$ be a nonnegative integer.
Let $(T,X)$ be a tree-decomposition of a graph $G$ of width strictly less than $w$.
Let $h = $\cc{$(2(a+2)w+2)b$}. If there is a cascade $\eta$ of height $h$
and size $s+k$ in $(T,X)$ such that $|\bigcap_{t \in V(T_{h})} X_{\eta(t)}| \geq k$,
then either there is a cascade $\eta'$ of height $a$
and size $s+k$ in $(T,X)$ such that
$|\bigcap_{t \in V(T_a)} X_{\eta'(t)}| \geq k+1$ or there is an injective cascade $\eta'$ of height $b$,
size $s+k$ and common intersection set of size $k$ in $(T,X)$.
\end{lemma}

\begin{proof}
We may assume that
\begin{itemize}
\item[$(*)$]
there does not exist a cascade $\eta'$ of height $a$ and size $s+k$ in $(T,X)$ such that
$|\bigcap_{t \in V(T_a)} X_{\eta'(t)}| \geq k+1$.
\end{itemize}
Let $F = \bigcap_{t \in V(T_{h})} X_{\eta(t)}$. By $(*)$, $|F| = k$.
We claim the following.
\begin{claim}
\label{cl:1}
For every vertex $t\in V(T_h)$ at height $l\le h-a-2$ and every $u\in X_{\eta(t)}-F$
  there exists a descendant $t'\in V(T_h)$ of $t$ at height at most
$l+a+2$ such that $u\not\in X_{\eta(t')}$.
\end{claim}

\noindent
To prove the claim let $u\in X_{\eta(t)}-F$.
By $(*)$ in the subtree of $T_{h}$ consisting of $t$ and its descendants
there is a vertex $t'$ of height at most $l+a + 2$ such that  $u\not\in X_{\eta(t')}$.
This proves the claim.
\bigskip

We use the previous claim to deduce the following generalization.

\begin{claim}
\label{cl:2}
For every vertex $t\in V(T_h)$ at height $l\le h-(a+2)w$ there exists a descendant $t'\in V(T)$ of $t$ at height
at most $l+(a+2)w$ such that $X_{\eta(t)}\cap X_{\eta(t')}=F$.
\end{claim}

\noindent
To prove the claim let $X_{\eta(t)}\backslash F = \{u_{1}, u_{2},\ldots, u_{p}\}$,
where $p\le w$. By Claim~\ref{cl:1}
 there exists a descendant $t_1\in V(T)$ of $t$ at height at most
$l+a+2$ such that $u_1\not\in X_{\eta(t')}$.
By another application of Claim~\ref{cl:1}
 there exists a descendant $t_2\in V(T)$ of $t_1$ at height at most
$l+2(a+2)$ such that $u_2\not\in X_{\eta(t')}$.
By (W2) $u_1\not\in X_{\eta(t')}$.
By continuing to argue in the same way we finally arrive at a vertex $t_p$
that is a descendant of $t$ at height at most $l+(a+2)p$ such that
$X_{\eta(t)}\cap X_{\eta(t_p)}=F$.
Thus $t_p$ is as desired.
This proves the claim.
\bigskip

Let  $x_0\in V(T_h)$ be the minor root of $T_h$.
By Claim~\ref{cl:2} and (W2) there exists
a major vertex $x\in V(T)$ at height at most $(a+2)w+1$ such that $X_{\eta(x_0)}\cap X_{\eta(x)}=F$.
Let   $y_1$ and $y_2$
be the children of $x$. By  Claim~\ref{cl:2} and (W2) there exists, for $i=1,2$, a minor vertex  $x_i\in V(T_h)$
at height  at most $2(a+2)w+2$  that  is a descendant of $y_i$ and such that $X_{\eta(x_i)}\cap X_{\eta(x)}=F$.
Let $r$ be the major root of $T_1$, and let $(t_0,t_1,t_2)$ be its trinity.
We define $\eta_1:T_1\emb T$ by $\eta_1(t_i)=\eta(x_i)$  for $i=0,1,2$ and $\eta_1(r)=\eta(x)$.
Then $\eta_1$ is an injective cascade of  height one and size $s+k$ in $(T,X)$ with common intersection set $F$.
If $b=1$, then $\eta_1$ is as desired, and so we may assume that \cc{$b>1$}.

%

\cc{Assume now that for some positive integer $l<b$ we have constructed an injective cascade $\eta_l:T_l\emb T$
of height $l$ and size $s+k$ with common intersection set $F$ in
$(T, X)$ such that for every leaf $t_0$ of $T_l$ other than the minor root
there exists a vertex $x_0\in V(T_h)$ at height $(2(a+2)w+2)l$ such that the image under $\eta_l$
of every vertex on the path in $T_l$ from the minor root to $t_0$ belongs to the path in $T$
from the image under $\eta$ of the root of $T_h$ to $\eta(x_0)$.
Our objective is to extend $\eta_l$ to an injective cascade
$\eta_{l+1}$ of height $l+1$, size $s+k$, and common intersection set $F$ in $(T, X)$
with the same property.
To that end let $\eta_{l+1}(t)=\eta_l(t)$ for all $t\in V(T_l)$, let $t_0$ be a leaf of $T_l$
other than the minor root, and let $x_0$ be as earlier in the paragraph.
By Claim~\ref{cl:2} and (W2) there exists a descendant $x$ of $x_0$
at height at most $(2(a+2)w+2)l + (a+2)w+1$ such that $x$ is major and $X_{\eta_l(t_0)}\cap X_{\eta(x)}=F$.
Let $y_1$ and $y_2$ be the children of $x$.
By  Claim~\ref{cl:2} and (W2) there exists, for $i=1,2$, a minor vertex  $x_i\in V(T_h)$
at height  at most $(2(a+2)w+2)(l+1)$  that  is a descendant of $y_i$ and such that $X_{\eta(x_i)}\cap X_{\eta(x)}=F$.
Let $r$ be the  child of $t_0$ in $T_{l+1}$, and let $(t_0,t_1,t_2)$ be its trinity.
We define  $\eta_{l+1}(t_i)=\eta(x_i)$  for $i=1,2$ and $\eta_{l+1}(r)=\eta(x)$.
This completes the definition of $\eta_{l+1}$.}

\cc{Now $\eta_b$ is as desired.}
\end{proof}

\begin{lemma}
\label{lemma4.8}
For any two positive integers $h$ and $w$,
there exists a positive integer $p = p(h, w)$ such that
if $G$ is a graph of tree-width less than $w$ and path-width at least $p$,
then in any tree-decomposition $(T, X)$ of $G$ that has width less than $w$ and
satisfies (W4), there is an injective cascade of height $h$.
\end{lemma}

\begin{proof}
%
\cc{Let $a_{w}= 0$, and for $k=w-1,\ldots,0$ let $a_k = (2(a_{k+1} + 2)w + 2)h$.
Let $p$ be the integer in Lemma~\ref{lemma4.6} for input integers $a_0$ and $w$.
We claim that $p$ satisfies the conclusion of the lemma.
To see that let $(T,X)$ be a tree-decomposition of $G$ of width less than~$w$ satisfying (W4).
By Lemma~\ref{lemma4.6}, there exists a cascade $\eta$ of height $a_0$ in $(T,X)$.
Let $k\in\{0,1,\ldots,w\}$ be the maximum integer such that
there exists a cascade $\eta':T_{a_{k}}\emb T$
satisfying $|\bigcap_{t \in V(T_{a_k})} X_{\eta'(t)}| \geq k$.
Such an integer exists, because $k=0$ satisfies those requirements
and $k<w$ because of (W4) and because the width of $(T,X)$ is less than $w$.
The maximality of $k$ implies that there does not exist a cascade $\eta'':T_{a_{k+1}}\emb T$
satisfying $|\bigcap_{t \in V(T_{a_{k+1}})} X_{\eta''(t)}| \geq k+1$.
Thus the lemma follows from Lemma~\ref{lemma4.7}.}
\end{proof}

\section{Ordered Cascades}
\label{sec:ordered}
The main result of this section, Theorem~\ref{thm:regular},
states that every $2$-connected graph of big path-width
and bounded tree-width admits a \td\ of bounded width and a cascade with linkages that are minimal.

Let $(T,X)$ be a tree-decomposition of a graph $G$, and let $\eta$ be
an injective cascade in $(T,X)$ with common intersection set $I$.
Assume the size of $\eta$ is $|I| + s$.
Then we say $\eta$ is {\em  ordered}  if for every minor vertex $t\in V(T_h)$
there exists a bijection $\xi_t:\{1,2,\ldots,s\}\to X_{\eta(t)}- I$
 such that for every major vertex $t_0$
with trinity $(t_1, t_2, t_3)$, there exist $s$ disjoint paths $P_1,P_2, \ldots,P_s$ in $G\backslash I$
such that the path $P_i$ has ends $\xi_{t_1}(i)$  and $\xi_{t_2}(i)$, and there exist $s$ disjoint paths
$Q_1,Q_2, \ldots,Q_s$ in $G\backslash I$ such that the path $Q_i$ has ends $\xi_{t_1}(i)$  and $\xi_{t_3}(i)$.
In that case we say that $\eta$  is an {\em  ordered cascade with orderings $\xi_t$}.
We say that the set of paths $P_1,P_2, \ldots,P_s$ is a {\em left  $t_0$-linkage with respect to $\eta$},
and that the set of paths $Q_1,Q_2, \ldots,Q_s$ is a {\em right  $t_0$-linkage with respect to $\eta$}.

 We will need to fix a left and a right
$t_0$-linkage for every major vertex $t_0\in V(T_h)$; when we do so we will indicate that by
saying that $\eta$ is an {\em ordered  cascade in $(T,X)$ with orderings $\xi_t$ and specified linkages},
and we will refer to the specified linkages as the
{\em left specified $t_0$-linkage} and the {\em right specified $t_0$-linkage}.
We will denote the  left specified $t_0$-linkage by $P_1(t_0),P_2(t_0),\ldots,P_s(t_0)$ and
the  right  specified $t_0$-linkage by $Q_1(t_0),Q_2(t_0),\ldots,Q_s(t_0)$.
We say that the specified $t_0$-linkages are {\em minimal} if
for
every set of disjoint paths
$P_1,P_2, \allowbreak\ldots,P_s$ in $G\backslash I$ from $X_{\eta(t_1)}-I$ to  $X_{\eta(t_2)}-I$ such that $\xi_{t_1}(i)$
is an end of $P_i$ (let the other end be $p_i$) and every set of disjoint paths $Q_1,Q_2, \ldots,Q_s$
in $G\backslash I$ from $X_{\eta(t_1)}-I$ to  $X_{\eta(t_3)}-I$  such that $\xi_{t_1}(i)$
is an end of $Q_i$  (let the other end be $q_i$)
we have
\begin{equation}
\label{c:min}
\left|E\left(\bigcup (x_iP_i p_i\cup x_iQ_i q_i)\right)\right|\ge
\left|E\left(\bigcup (y_iP_i(t_0)\xi_{t_2}(i)\cup y_iQ_i(t_0)\xi_{t_3}(i))\right)\right|,
\end{equation}
where the unions are taken over $i\in \{1, 2, ..., s\}$,
$x_i$ is the first vertex from $\xi_{t_1}(i)$ that
$P_i$ departs from $Q_i$, and $y_i$ is the first vertex
from $\xi_{t_1}(i)$ that $P_i(t_0)$ departs from
$Q_i(t_0)$.

\begin{lemma}
\label{lemma4.9}
Let $h$ and $s$ be two positive integers, and let  $\eta:T_h\emb T$ be an injective cascade
of height $h$ and size $s$ in a \cd{linked} \td\ $(T, X)$ of a graph $G$.
Then the cascade $\eta$ can be turned into an ordered cascade with specified $t_0$-linkages that are
minimal for every major vertex $t_0\in V(T_h)$.
\end{lemma}

\begin{proof}
Let $s':=s-|I|$.
To show that $\eta$ can be made ordered let $r$ be the minor root of $T_h$, let
$\xi_r:\{1,2,\ldots,s'\}\to X_{\eta(r)}-I$ be arbitrary, assume that for some integer $l\in\{0,1,\ldots,h-1\}$
 we have already
constructed $\xi_t:\{1,2,\ldots,s'\}\to X_{\eta(t)}-I$ for all minor vertices $t\in V(T_h)$ at height
at most $l$, let  $t\in V(T_h)$ be a minor vertex at height exactly $l$, let $t_0$ be its child,
and let $(t,t_1,t_2)$ be the trinity at $t_0$.
By condition (W3) there exist $s'$ disjoint paths $P_1,P_2,\ldots,P_{s'}$ in $G\backslash I$ from $X_{\eta(t)}-I$
to $X_{\eta(t_1)}-I$ and  $s'$ disjoint paths $Q_1,Q_2,\ldots,Q_{s'}$ in $G\backslash I$ from $X_{\eta(t)}-I$
to $X_{\eta(t_2)}-I$.
We may assume that $\xi_t(i)$ is an end of $P_i$ and $Q_i$  and we define $\xi_{t_1}(i)$ and $\xi_{t_2}(i)$ to
be their other ends, respectively. We may also assume that these paths satisfy the minimality condition~(\ref{c:min}).
It follows that $\eta$ is an ordered cascade with orderings $\xi_t$ and  specified $t_0$-linkages that are
minimal for every major vertex $t_0\in V(T_h)$.
\end{proof}

%

Let $h,h'$ be integers. We say that a \he\ $\gamma:T_{h'}\emb T_h$
is {\em monotone} if
\begin{itemize}
\item $t$ is a major vertex of $T_{h'}$ with trinity $(t_1,t_2,t_3)$, then $\gamma(t_2)$ is the
left neighbor of $\gamma(t)$ and $\gamma(t_3)$ is the
right neighbor of $\gamma(t)$,
and
\item the image under $\gamma$ of the minor root of $T_{h'}$ is the minor root of $T_h$.
\end{itemize}

\begin{lemma}
\label{lemma4.10}
For every two integers $a \ge 1$ and $k \ge 1$
there exists an integer $h = h(a, k)$ such that the following holds.
Color the major vertices of $T_h$ using $k$ colors.
Then there exists a monotone \he\ $\eta: T_{a} \emb T_h$ such that the major vertices of
$T_{a}$ map to major vertices of the same color in $T_h$.
\end{lemma}

\begin{proof}
Let $c$ be one of the colors. We will prove by induction on $k$ and subject to that
by induction on $b$ that
there is a function $h=g(a,b,k)$ such that there is either  a monotone \he\ $\eta: T_{a} \emb T_h$
 such that the major vertices of
$T_{a}$ map to major vertices of the same color in $T_h$, or  a monotone \he\ $\eta: T_{b} \emb T_h$
such that the major vertices of
$T_{b}$ map to major vertices of  color  $c$ in $T_h$.
In fact, we will show that $g(a,b,1)=a$, $g(a,1,k+1)\le g(a,a,k)$ and $g(a,b+1,k+1)\le g(a,b,k+1)+g(a,a,k)$.

The assertion holds for $k = 1$ by letting $h=a$ and letting $\eta$ be the identity mapping.
Assume the statement is true for some $k \ge 1$,  let
the major vertices of $T_h$ be colored using
$k+1$ colors, and let $c$ be one of the colors.
If $b=1$, then if $T_h$ has a major vertex colored $c$, then the second alternative holds;
otherwise at most $k$ colors are used and the assertion  follows by induction on $k$.

We may therefore assume that the assertion holds for some integer $b\ge1$ and we
must prove it for $b+1$.
To that end we may assume that $T_h$ has a major vertex $t_0$ colored $c$ at height
at most $g(a,a,k)$, for otherwise the assertion follows by induction on $k$.
Let the trinity at $t_0$ be $(t_1,t_2,t_3)$.
For $i=2,3$ let $R_i$ be the subtree of  $T_h$ with minor root~$t_i$.
If  for some $i\in\{2,3\}$ there exists  a monotone \he\ $ T_{a} \emb R_i$
 such that the major vertices of
$T_{a}$ map to major vertices of the same color in $T_h$, then the statement holds.
We may therefore assume that for $i\in\{2,3\}$ there exists  a monotone \he\ $\eta_i: T^i_{b} \emb R_i$
 such that the major vertices of
$T^i_{b}$ map to major vertices of color $c$, the major root of $T_{b+1}$ is $r_0$,
the trinity at $r_0$ is $(r_1,r_2,r_3)$ and $T^i_b$ is the subtree of $T_{b+1}-\{r_0,r_1\}$
with minor root $r_i$.
Let $\eta:T_{b+1}\emb T_h$ be defined by $\eta(t)=\eta_i(t)$ for $t\in V(T^i_b)$,
$\eta(r_0)=t_0$ and $\eta(r_1)$ is defined to be the minor root of $T_h$.
Then $\eta:T_{b+1}\emb T_h$ is as desired.
This proves the existence of the function $g(a,b,k)$.

Now $h(a,k)=g(a,a,k)$ is as desired.
%
%
\end{proof}

Let $G$ be a graph, let $v\in V(G)$ and for $i=1,2,3$ let $P_i$
be a path in $G$ with ends $v$ and $v_i$
such  that the paths $P_1,P_2,P_3$ are pairwise  disjoint, except for $v$.
Assume that at least two of the paths $P_i$ have length at least one.
We say that $P_1\cup P_2\cup P_3$ is a {\em tripod} with
{\em center} $v$ and {\em feet} $v_1,v_2,v_3$.

Let $(T,X)$ be a tree-decomposition of a graph $G$, and let $\eta:T_h\emb T$ be
an injective cascade in $(T,X)$ with common intersection set $I$.
Let $t_0\in V(T_h)$ be
a major vertex, and let $(t_1,t_2,t_3)$ be the trinity at $t_0$.
We define the {\em $\eta$-torso at $t_0$} as the subgraph of $G$ induced by $\bigcup X_t- I$, where the union is taken
 over all $t$ in $V(T)$ such that the unique path in $T$ from $t$ to $\eta(t_0)$ does not contain
$\eta(t_1)$,$\eta(t_2)$, or $\eta(t_3)$ as an internal vertex.

Let $s > 0$ be an integer. Let $(T,X)$ be  a tree-decomposition of a graph $G$,
let $\eta:T_h\emb T$ be an ordered cascade in $(T,X)$ of size $|I| + s$ and with
 orderings $\xi_t$, where $I$ is
the common intersection set of $\eta$. Let $t_0\in V(T_h)$  be a major vertex,
let $(t_1,t_2,t_3)$ be the trinity at $t_0$, let $G'$ be the $\eta$-torso at $t_0$, and let
$i,j\in\{1,2,\ldots,s\}$ be distinct.
We say that {\em $t_0$ has property $A_{ij}$ in $\eta$} if there exist disjoint tripods $L_i, L_j$ in $G'$
such that for each $m\in\{i,j\}$ the tripod $L_m$ has feet $\xi_{t_1}(m),\xi_{t_2}(m_2),\xi_{t_3}(m_3)$
for some $m_2, m_3 \in \{i,j\}$.

We say that {\em $t_0$ has property $B_{ij}$ in $\eta$} if there exist
vertices $v_{x,y}$ for all $x\in \{i,j\}, y\in \{1,2,3\}$,
and tripods $L_i, L_j$ in $G'$ with centers $c_i, c_j$ such that
\begin{itemize}
\item for each $y\in \{1,2,3\}$, $\{v_{i,y}, v_{j,y}\} = \{\xi_{t_y}(i),\xi_{t_y}(j)\}$
\item for each $m\in\{i,j\}$, $L_m$ has feet $v_{m,1}, v_{m,2}, v_{m,3}$
\item $L_{i} \cap L_{j} = c_{i}L_{i}v_{i,3} \cap c_{j}L_{j}v_{j,2}$ and it is a path that does not
contain $c_{i}, c_{j}$.
\end{itemize}

We say that {\em $t_0$ has property \cd{$C_{ij}$} in $\eta$} if there exist three pairwise
disjoint paths $R_i,R_j, R_{ij}$ and a path $R$ in $G'$ such that
the ends of $R_i$ are $\xi_{t_1}(i)$ and $\xi_{t_2}(i)$,
the ends of $R_j$ are $\xi_{t_1}(j)$ and $\xi_{t_3}(j)$,
the ends of $R_{ij}$ are $\xi_{t_2}(j)$ and $\xi_{t_3}(i)$, and
$R$ is internally disjoint from $R_i, R_j, R_{ij}$ and connects two of these three paths.
We will denote these paths as $R_i(t_0), R_j(t_0), R_{ij}(t_0), R(t_0)$ when we want to
emphasize they are in the torso at the major vertex $t_0$.

We say that the path $P_i$ of a left or  right  $t_0$-linkage is {\em confined} if it is a subgraph of the $\eta$-torso at $t_0$.

Now let $\eta: T_h \emb T$ be an ordered cascade in
$(T,X)$ with orderings $\xi_t$ and specified linkages.
Let $t_0\in V(T_h)$  be a major vertex with trinity $(t_1,t_2,t_3)$,
and let $P_1,P_2, \ldots,P_s$ be the left specified $t_0$-linkage.
We define $A_{t_0}$  to be the set of integers $i\in\{1,2,\ldots,s\}$
such that the path $P_i$ is confined,
and we define $B_{t_0}$ in the same way but using the
right specified $t_0$-linkage instead.
Define $C_{t_0}$ as the set of all triples $(i, l, m)$
such that $i\in\{1,2,\ldots,s\}$, the path $P_i$
is not confined and  when  following $P_i$ from $\xi_{t_1}(i)$,
it exits the $\eta$-torso at $t_0$ for the first time at
$\xi_{t_3}(l)$ and re-enters the $\eta$-torso at $t_0$
for the last time at $\xi_{t_3}(m)$.
Let $D_{t_0}$ be defined similarly, but using the  right $t_0$-linkage instead.
We call the sets $A_{t_0},B_{t_0}$, $C_{t_0}$ and
$D_{t_0}$ the {\em confinement sets for $\eta$ at $t_0$
with respect to the specified linkages.}

Let $A_{t_0}$ and $B_{t_0}$ be the confinement sets for $\eta$ at $t_0$.
We say that {\em $t_0$ has property \cd{$C$} in $\eta$} if $s$ is even,
$A_{t_0}$  and $B_{t_0}$ are disjoint
and both have size $s/2$, and there exist disjoint paths
$R_1,R_2,\ldots,R_{3s/2}$ in $G'$
in such a way that

\begin{itemize}
\item each $R_i$ is a subpath of both the left specified $t_0$-linkage and the right specified $t_0$-linkage,
\item for $i\in A_{t_0}$, the path $R_i$ has ends $\xi_{t_1}(i)$ and $\xi_{t_2}(i)$,
\item for $i\in B_{t_0}$ the path $R_i$ has ends $\xi_{t_1}(i)$ and $\xi_{t_3}(i)$, and
\item for $i=s+1,s+2,\ldots,3s/2$ the path $R_i$ has one end $\xi_{t_2}(k)$ and the other and $\xi_{t_3}(l)$
for some $k\in B_{t_0}$ and $l\in A_{t_0}$.
\end{itemize}

Let $(T,X)$ be  a tree-decomposition of a graph $G$,
let $\eta:T_h\emb T$ be
a cascade in $(T,X)$  and let $\gamma:T_{h'}\emb T_h$  be a monotone \he.
Then the composite mapping $\eta':=\eta\circ\gamma: T_{h'} \emb T$ is a cascade in $(T,X)$ of height $h'$,
and we will call it a {\em subcascade of $\eta$}.

\begin{lemma}
\label{lem:subcasc}
Let $(T,X)$ be  a tree-decomposition of a graph $G$,
let $\eta:T_h\emb T$ be
an ordered cascade in $(T,X)$ with orderings $\xi_t$, specified linkages and common intersection set $I$,
let $\gamma:T_{h'}\emb T_h$  be a monotone \he,
and let $\eta':=\eta\circ\gamma: T_{h'} \emb T$ be a subcascade of $\eta$ of height $h'$.
Then for every major vertex  $t_0\in V(T_{h'})$
\begin{enumerate}
\item[\rm(i)]
 $\eta'$ is an ordered cascade with orderings $\xi_{\gamma(t)}$ and common intersection set~$I$,
\item[\rm(ii)]
 if the vertex $\gamma(t_0)$ has property $A_{ij}$ ($B_{ij}$, $C_{ij}$, resp.) in $\eta$,
then $t_0$  has property $A_{ij}$ ($B_{ij}$, $C_{ij}$, resp.) in $\eta'$.
\end{enumerate}
Furthermore, 
the specified linkages for $\eta'$
may be chosen in such a way that
\begin{enumerate}
\item[\rm(iii)]
$(A_{t_0}, B_{t_0}, C_{t_0}, D_{t_0}) =
(A_{\gamma(t_0)}, B_{\gamma(t_0)}, C_{\gamma(t_0)}, D_{\gamma(t_0)})$,
\item[\rm(iv)]
the vertex $t_0$ has property C in $\eta'$
if and only if $\gamma(t_0)$ has property C in $\eta$, and
\item[\rm(v)]
 if the specified linkages for $\eta$ are minimal, then the specified linkages for $\eta'$
are minimal.
\end{enumerate}
\end{lemma}

\begin{proof}
For each major vertex $t\in V(T_{h'})$ or $t\in V(T_h)$ we denote its trinity by $(t_1(t), t_2(t),\allowbreak t_3(t))$.
Assume $t_0$ is a major vertex of $T_{h'}$. Let $v_0 = \gamma(t_1(t_0)), v_1, \ldots, v_k = t_1(\gamma(t_0))$
be the minor vertices on $T_h[v_0, v_k]$. Let $U$ be the union of the left (or right) linkage from $X_{\eta(v_i)}- I$
to $X_{\eta(v_{i+1})}- I$ for all $i \in \{0, 1,\ldots, k-1\}$ depending on whether $v_{i+1}$ is a left (or right)
neighbor of its parent. Let $P$ be the left
specified $\gamma(t_0)$-linkage and $Q$ be the right specified $\gamma(t_0)$-linkage.
Then $U\cup P$ is a left $t_0$-linkage and $U\cup Q$ is a right $t_0$-linkage.
We designate  $U\cup P$ to be the left specified $t_0$-linkage and  $U\cup Q$ to be the right specified $t_0$-linkage.
It is easy to see that this choice satisfies the conclusion of the lemma.
\end{proof}

Let $(T,X)$ be  a tree-decomposition of a graph $G$,
and let $\eta$ be an ordered cascade with specified linkages
in $(T,X)$ of height $h$ and size $|I| + s$,
where $I$ is the common intersection set.
We say that $\eta$ is {\em regular} if
there exist sets $A,B\subseteq\{1,2,\ldots,s\}$, and sets $C$ and $D$ such that
the confinement sets $A_{t_0}$, $B_{t_0}$, $C_{t_0}$
and $D_{t_0}$ satisfy $A_{t_0}=A$, $B_{t_0}=B$, $C_{t_0}=C$
and $D_{t_0}=D$ for
every  major vertex $t_0 \in V(T_{h})$.

\begin{lemma}
\label{lemma4.11}
For every two positive integers $a$ and $s$ there exists a positive integer
$h = h(a, s)$ such that the following holds.
Let $(T,X)$ be a \cd{linked} tree-decomposition of a graph $G$.
If there exists an injective cascade $\eta$ of height $h$  in $(T, X)$,
then there exists a regular cascade $\eta':T_a\emb T$ of height $a$ in $(T, X)$
with specified $t_0$-linkages that are minimal for every major vertex $t_0\in V(T_a)$
such that $\eta'$ has the same size and common intersection set as $\eta$.
\end{lemma}

\begin{proof}
Let $\eta$ be an injective cascade  of size $|I| + s$ and height $h$ in $(T, X)$,
where we will specify $h$ in a moment.
By Lemma~\ref{lemma4.9} $\eta$ can be turned into an ordered cascade with specified $t_0$-linkages that are
minimal for every major vertex $t_0\in V(T_h)$.
For every major vertex $t_0 \in V(T_h)$, the number of possible quadruples
$(A_{t_0}, B_{t_0}, C_{t_0}, D_{t_0})$ is a finite number $k = k(s)$ that depends only on $s$.

Consider each choice of $(A_{t_0}, B_{t_0}, C_{t_0}, D_{t_0})$ as a color;
then by Lemma~\ref{lemma4.10},
there exists a positive integer $h = h(a, k)$ such that there exists a monotone homeomorphic
embedding $\gamma: T_a \emb T_h$ such that the quadruple
$(A_{\gamma(t)}, B_{\gamma(t)}, C_{\gamma(t)}, D_{\gamma(t)})$
for $\eta$ is the same for every $t \in V(T_a)$.
%
Now, let $\eta' = \eta \circ \gamma : T_a \rightarrow T$. Then
$\eta'$ is as desired by Lemma~\ref{lem:subcasc}.
%
%
\end{proof}


The following is the main result of this section.

\begin{theorem}
\label{thm:regular}
For any two positive integers $a$ and $w$, there exists a positive integer $p = p(a, w)$
such that the following holds.
Let $G$ be a $2$-connected graph of tree-width less than $w$ and path-width at least $p$.
Then $G$ has a tree-decomposition $(T,X)$ such that:
\begin{itemize}
\item $(T, X)$ has width less than $w$,
\item $(T, X)$ satisfies (W1)--(W7), and
\item for some $s$, where $2 \leq s \leq w$, there exists a regular cascade $\eta:T_a\emb T$
of height $a$ and size $s$ in $(T, X)$ with specified $t_0$-linkages that are
minimal for every major vertex $t_0\in V(T_a)$.
\end{itemize}
\end{theorem}

\begin{proof}
Given  positive integers $a$ and $w$  let $h$ be as in Lemma~\ref{lemma4.11},
and let $p = p(h, w)$ be as in  Lemma~\ref{lemma4.8}. We claim that $p$
satisfies the conclusion of the theorem. To see that
let $G$ be a graph of tree-width less than $w$ and path-width at least $p$.
By Theorem~\ref{w7}, $G$ admits a tree-decomposition $(T, X)$ of width less than $w$
 satisfying (W1)--(W7). By Lemma~\ref{lemma4.8} there is
an injective cascade
of height $h$ in $(T,X)$.
Let $s$ be the size of this cascade, then
$s \leq w$. If $G$ is $2$-connected, then $s \geq 2$.
The last conclusion of the theorem follows from Lemma~\ref{lemma4.11}.
\end{proof}

\section{Taming Linkages}
\label{sec:taming}

Lemma~\ref{lemma4.15}, the main result of this section, states that there are essentially only two 
types of linkage.

Let $s > 0$ be an integer. Let $(T,X)$ be  a tree-decomposition of a graph $G$,
let $\eta:T_h\emb T$ be an ordered cascade in $(T,X)$ of size $|I| + s$ and with orderings $\xi_t$, where $I$ is
the common intersection set of $\eta$. Let $t_0\in V(T_h)$  be a major vertex,
let $(t_1,t_2,t_3)$ be the trinity at $t_0$, let $G'$ be the $\eta$-torso at $t_0$, and let
$i,j\in\{1,2,\ldots,s\}$ be distinct.
We say that {\em $t_0$ has property $AB_{ij}$ in $\eta$} if there exist
disjoint paths $L_i,L_j$ and disjoint paths $R_i, R_j$ in $G'$
such that the two ends of $L_m$ are $\xi_{t_1}(m)$ and $\xi_{t_2}(m)$ for each
$m\in\{i,j\}$ and the two ends of $R_m$ are $\xi_{t_1}(m)$ and $\xi_{t_3}(m)$ for each $m\in\{i,j\}$.

If $P$ is a path and $u,v \in V(P)$, then by $uPv$ we denote the subpath of $P$
with ends $u$ and $v$.

\begin{lemma}
\label{lemma4.13}
Let $(T,X)$ be  a tree-decomposition of a graph $G$. Let $\eta:T_1\emb T$ be an
ordered cascade in $(T,X)$ with orderings $\xi_t$ of height one and size $s+|I|$,
where $I$ is the common intersection set. Let $t_0$ be the major vertex in $T_1$,
and let $i,j\in\{1,2,\ldots,s\}$ be distinct.
If $t_0$ has property $AB_{ij}$ in $\eta$, then $t_0$ has either property $A_{ij}$
or property $B_{ij}$ in $\eta$.
\end{lemma}

\begin{proof}
Let $(t_1,t_2,t_3)$ be the trinity at $t_0$. Let $G'$ be the $\eta$-torso at $t_0$.
Since $t_0$ has property $AB_{ij}$ in $\eta$, there exist disjoint paths $L_i,L_j$ and
disjoint paths $R_i, R_j$ in $G'$ such that two endpoints of $L_m$ are $\xi_{t_1}(m)$
and $\xi_{t_2}(m)$ for all $m\in\{i,j\}$, and two endpoints of $R_m$ are $\xi_{t_1}(m)$
and $\xi_{t_3}(m)$ for all $m\in\{i,j\}$.

We may choose $L_i,L_j, R_i, R_j$ such that
$|E(L_i) \cup E(L_j) \cup E(R_i) \cup E(R_j)|$ is as small as possible.

Let $x_k = \xi_{t_1}(k)$
and $z_k = \xi_{t_3}(k)$ for  $k \in \{i, j\}$.
Starting from $z_i$, let $a$ be the first vertex where $R_i$ meets $L_i \cup L_j$,
and starting from $z_j$, let $b$ be the first vertex where $R_j$ meets $L_i \cup L_j$.
If $a$ and $b$ are not on the same path (one on $L_i$ and the other on $L_j$),
then
by considering $L_i, L_j$ and the parts of $R_i$ and $R_j$ from $z_i$ to $a$
and from $z_j$ to $b$ we see that $t_0$ has property $A_{ij}$ in $\eta$.

If $a$ and $b$ are on the same path, then we may assume they are on $L_i$.
We may also assume that $a\in L_i[y_i,b]$.
Then following $R_i$ from $a$ away from $z_i$, the paths $R_i$ and $L_i$ eventually split; let $c$
be the vertex where the split occurs.
In other words, $c$ is such that $aL_ic\cap aR_ic$ is  a path and its length is maximum.
Let $d$  be the first vertex on $cR_ix_i\cap(L_i\cup L_j)-\{c\}$ when traveling on $R_i$ from $c$ to $x_i$.
If $d\in V(L_i)$, then by replacing $cL_id$ by $cR_id$ we obtain a contradiction to the choice of
 $L_i,L_j, R_i, R_j$.
Thus $d\in V(L_j)$.
Now $L_i,L_j$ and the paths $z_iR_id$ and $z_jR_jb$ show that
$t_0$ has property $B_{ij}$ in $\eta$.
%
\end{proof}

Let $(T,X)$ be a tree-decomposition of a graph $G$
and let  $\eta:T_h\emb T$ be an injective cascade
in $(T,X)$
 of height $h$ and size $|I| + s$, where $I$ is the common intersection set.
Let $v$ be a vertex of $T_h$ and let $Y$  consist of $\eta(v)$ and the  vertex-sets of all components of
$T\backslash\eta(v)$ that  do not contain the image under $\eta$ of the minor root of  $T_h$.
Let $H$ be the subgraph of $G$  induced by $\bigcup_{t\in Y} X_t-I$.
We will call $H$ the {\em outer graph at $v$}.

\begin{lemma}
\label{lem:W6}
Let $(T,X)$ be a tree-decomposition satisfying (W6) of a graph $G$
and let  $\eta:T_h\emb T$ be an \cd{ordered} cascade
in $(T,X)$
of height $h$ and size $|I| + s$, where $I$ is the common intersection set.
Let $v$ be a minor vertex of $T_h$ at  height at most $h-1$,  let $H$ be the outer graph at $v$, and let $x,y\in X_{\eta(v)}$.
Then there exists a  path  of length at least two with ends $x$ and $y$ and every internal vertex in $V(H)-X_{\eta(v)}$.
\end{lemma}

\begin{proof}
Let $v_0$ be the child of $v$, let $v_1$ be a child of $v_0$, and let $B$ be the component of $T -\eta( v)$ that contains $\eta(v_1)$. We show that $x$ is $B$-tied.
This is obvious if $x\in I$, and so we  may assume that $x\not\in I$.
Since $\eta$ is ordered, there exist $s$ disjoint paths from $X_{\eta(v)}-I$  to
$X_{\eta(v_1)} - I$ in $G\backslash I$. It follows that each of the paths uses exactly one vertex of $X_{\eta(v)}-I$,
and that vertex is its end.
Let $P$ be the one of those paths that ends in $x$, and let $x'$  be the neighbor of $x$ in $P$.
The vertex $x'$ exists, because $X_{\eta(v)}\cap X_{\eta(v_1)}=I$.
By (W1) there exists a vertex $t\in V(T)$ such that $x,x'\in X_t$. Since $P-x$ is disjoint from $X_{\eta(v)}$, it follows
from Lemma~\ref{cutbag} applied to the path $P-x$ and vertices $t$ and $\eta(v_1)$ of $T$ that $t\in V(B)$.
Thus $x$ is $B$-tied and the same argument shows that so is $y$. Hence the lemma follows from (W6).
\end{proof}

We will refer
to a path as in Lemma~\ref{lem:W6} as a W6-{\em path}.

Let $h,h'$ be integers. We say that a \he\ $\gamma:T_{h'}\emb T_h$
is {\em weakly monotone} if  for every two vertices
$t,t'\in V(T_{h'})$
\begin{itemize}
\item if $t'$   is a descendant of $t$ in $T_{h'}$,
then the vertex $\gamma(t')$  is a descendant of $\gamma(t)$ in $T_{h}$
\item if $t$  is a minor vertex of $T_{h'}$,
then the vertex $\gamma(t)$  is minor in $T_{h}$.
\end{itemize}
Let $(T,X)$ be  a tree-decomposition of a graph $G$,
let $\eta:T_h\emb T$ be
a cascade in $(T,X)$  and let $\gamma:T_{h'}\emb T_h$  be a weakly monotone \he.
Then the composite mapping $\eta':=\eta\circ\gamma: T_{h'} \emb T$ is a cascade in $(T,X)$ of height $h'$,
and we will call it a {\em weak subcascade of $\eta$}.

\begin{lemma}
\label{lemma4.14}
Let $s \geq 2$ be an integer,
let $(T,X)$ be a tree-decomposition of a graph $G$ satisfying (W6),
and let \cb{$\eta:T_5\emb T$} be a regular cascade
in $(T,X)$ of height five and size $|I| + s$
with specified linkages that are minimal,
where $I$ is the common intersection set of $\eta$.
Then either there exists a weak subcascade $\eta':T_{1}\emb T$ of $\eta$ of height one
such that the unique major vertex of
$T_1$ has property
$A_{ij}$ or $B_{ij}$  in $\eta'$
for some distinct integers $i,j\in\{1,2,\ldots,s\}$,
or
the major root of $T_5$ has property $C$ in $\eta$.
\end{lemma}

\begin{proof}
We  will either construct a weakly monotone \he\ $\gamma: T_{1}\emb T_5$ such that in
$\eta' = \eta \circ \gamma$ the major root of $T_1$ will have property $AB_{ij}$ for some distinct $i, j \in \{1, 2, ..., s\}$,
or establish that the major root of $T_5$ has property $C$ in $\eta$.
By Lemma~\ref{lemma4.13} this will suffice.

Since $\eta$ is regular, there exist sets $A,B,C,D$  as in the definition of regular cascade.
Let $t_0$ be the unique major vertex of $T_1$ and let $(t_1, t_2, t_3)$  be its trinity.
Let $u_0$ be the major root of $T_5$ and let $(v_1, v_2, v_3)$ be its  trinity.
Let $u_1$ be the major vertex of $T_5$ of height one  that is adjacent to $v_3$ and let $(v_3, v_4, v_5)$ be its trinity.
Let us recall that for a major vertex $u$ of $T_5$ we denote the paths in the
 specified left $u$-linkage by  $P_i(u)$ and the paths in the specified right $u$-linkage by  $Q_i(u)$.
If there exist two distinct integers $i,j\in A\cap B$, then the paths $P_i(u_0), P_j(u_0),Q_i(u_0),Q_j(u_0)$
show that $u_0$ has property $AB_{ij}$ in $\eta$.
Let $\gamma: T_{1}\emb T_5$ be the \he\ that
maps $t_0,t_1,t_2,t_3$ to $u_0,v_1,v_2,v_3$, respectively.
Then  $\eta' = \eta \circ \gamma$ is as desired.
We may therefore assume that $|A\cap B|\le1$.


For $i\in\{1,2,\ldots,s\}-A$  the path $P_i(u_0)$ exits and re-enters the $\eta$-torso at $u_0$,
and it does so through  two distinct vertices of $X_{\eta(v_3)}$. But $|X_{\eta(v_3)}- I|=s$,
and hence $|A|\ge s/2$. Similarly $|B|\ge s/2$.
By symmetry we may assume that $|B|\ge|A|$. It follows that $|A|=\lceil s/2\rceil$, and hence for
$i\in\{1,2,\ldots,s\}-A$ and every major vertex $w$ of $T_5$ the path $P_i(w)$
exits and re-enters the $\eta$-torso at $w$ exactly once.
The set $C$ includes an element of the form $(i,l,m)$, which means that the vertices
$\xi_{w_1}(i),\xi_{w_3}(l),\xi_{w_3}(m),\xi_{w_2}(i)$ appear on the path $P_i(w)$ in the order listed.
Let $l_i:=l,m_i:=m$, $x_i(w):=\xi_{w_3}(l)$, $y_i(w):=\xi_{w_3}(m)$, $X_i(w):=\xi_{w_1}(i)P_i(w)x_{i}(w)$ and
\redsout{$Y_i(w):=x_{i}(w)P_i(w)\xi_{w_2}(i)$}%
\cb{$Y_i(w):=y_{i}(w)P_i(w)\xi_{w_2}(i)$}.
Thus $X_i(w)$ and $Y_i(w)$ are  subpaths of the $\eta$-torso at $w$.
We distinguish two main cases.

\medskip

\noindent
{\bf Main case 1:} $|A \cap B| = 1$.
Let $j$  be the unique element of $A \cap B$.
We claim that $B-A\ne\emptyset$. To prove the claim suppose for a contradiction that $B\subseteq A$.
Thus $|B|=1$, and since $|B|\ge|A|$ we have $|A|=1$, and hence $s=2$.
We may assume,  for the duration of this paragraph, that $A=B=\{1\}$.
The paths $P_1(u_0),X_2(u_0),Y_2(u_0)$ are pairwise disjoint, because they are  subgraphs of the specified left $u_0$-linkage.
 The path $Q_2(u_0)$ is unconfined, and  hence it has a subpath $R$  joining $\xi_{v_2}(1)$ and $\xi_{v_2}(2)$
in the outer graph at $v_2$. It follows that $P_1(u_0)\cup R\cup Y_2(u_0)$ and $X_2(u_0)$ are disjoint paths from $X_{\eta(v_1)}$
to $X_{\eta(v_3)}$, and it follows from the minimality of the specified $u_0$-linkage
that they form the specified right $u_0$-linkage,
 contrary to $1\in A$.
This proves the claim that $B-A\ne\emptyset$, and so we may select an element $i\in B-A$.

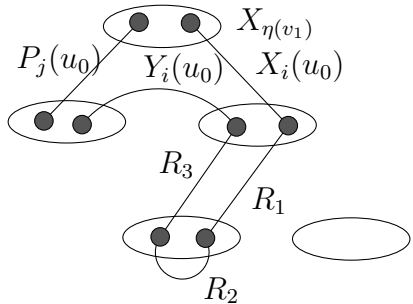
\begin{figure}[htb]
 \centering
\input{image10.tikz}
 \caption{First case of the construction of the path $R$.}
\label{fig:abcd1}
\end{figure}

Let us assume as a case that either $l_i\in A$ or $l_i\not\in B$.
In this case we let $\gamma$ map $t_0,t_1,t_2,t_3$ to $u_0,v_1,v_2,v_5$, respectively,
and we will prove that $t_0$ has property $AB_{ij}$ in $\eta'$.
To that end we need to construct two pairs of disjoint paths.
The first pair is $Q_i(u_0)\cup Q_i(u_1)$
and $Q_j(u_0)\cup Q_j(u_1)$.
The second pair will consist of $P_j(u_0)$ and another path \cb{from $\xi_{v_1}(i)$ to $\xi_{v_2}(i)$} which
\cb{is a subgraph of a walk that}
we are about to construct.
It will consist of $X_i(u_0)\cup Y_i(u_0)$ and a \redsout{path}%
\cb{walk} $R$
in the outer graph of $v_3$ with ends $x_i(u_0)$ and $y_i(u_0)$.
To construct the \redsout{path }%
\cb{walk} $R$ we will construct \redsout{three}%
paths
$R_1,R_2$ \cb{and a walk} $R_3$,  whose union will contain the desired \redsout{path}%
\cb{walk} $R$.
If $l_i\in A$, then we let $R_1:=P_{l_i}(u_1)$. If $l_i\not\in B$, then the path $Q_{l_i}(u_1)$ is unconfined,
and hence includes a subpath $R_1$ from $x_i(u_0)$  to $X_{\eta(v_4)}$  that is a subgraph of the $\eta$-torso at $u_1$.
We need to distinguish two  subcases  depending on whether $m_i\in B$.
Assume first that $m_i\not\in B$ and refer to Figure~\ref{fig:abcd1}. Then similarly as above the path $Q_{m_i}(u_1)$ is unconfined,
and hence includes a subpath $R_3$ from $y_i(u_0)$ to $X_{\eta(v_4)}$  that is a subgraph of the $\eta$-torso at $u_1$,
and we let $R_2$ be a W6-path in the outer graph at $v_4$ joining the ends of $R_1$ and $R_3$ in $X_{\eta(v_4)}$.
This completes the subcase $m_i\not\in B$, and so we may assume that $m_i\in B$.
In this  subcase we define $R_3:=Y_{i}(u_1)\cup Q_{m_i}(u_1)$ and we define $R_2$ as above.
See Figure~\ref{fig:abcd2}.
This completes the case that either $l_i\in A$ or  $l_i\not\in B$.

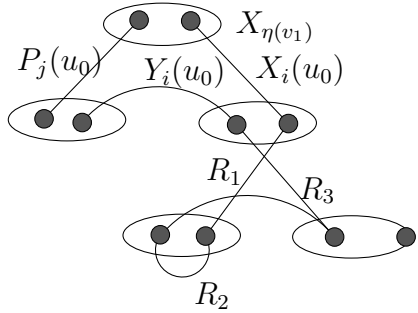
\begin{figure}[htb]
 \centering
\input{image11.tikz}
 \caption{Second case of the construction of the path $R$.}
\label{fig:abcd2}
\end{figure}

Next we consider the case $l_i\in B$ and $m_i\not\in A-B$.
We proceed similarly as in the previous paragraph, but with these exceptions:
the \he\ $\gamma$ will map $t_3$ to $v_4$, rather than $v_5$,
the first pair of disjoint paths will now be  $Q_i(u_0)\cup P_i(u_1)$ and $Q_j(u_0)\cup P_j(u_1)$,
and for the second pair we define
$R_1=Q_{l_i}(u_1)$,  $R_3=X_{m_i}(u_1)$ if $m_i\not\in A$ and  $R_3=Q_{m_i}(u_1)$ if $m_i\in B$,
and $R_2$ will be a W6-path in the outer graph
of $v_5$  joining the ends of $R_1$ and $R_3$.

Therefore assume that $l_i\in B-A$ and $m_i\in A-B$ for every $i\in B-A$.
Let $u_2$ be the major vertex of $T_5$ at height two whose trinity includes $v_5$
and assume its trinity is $(v_5, v_6, v_7)$.
\cb{Let $u_3$ be the major vertex of $T_5$ at height three whose trinity includes $v_7$
and assume its trinity is $(v_7, v_8, v_9)$.}
Let $\gamma$ map $t_0,t_1,t_2,t_3$ to $u_0,v_1,v_2,\cb{v_8}$, respectively.
Then $t_0$ also has property $AB_{ij}$ in $\eta'$. To see that
the first pair of disjoint paths is $Q_i(u_0)\cup Q_i(u_1)\cb{\cup Q_i(u_2)} \cup P_i(\cb{u_3})$
and $Q_j(u_0)\cup Q_j(u_1)\cb{\cup Q_j(u_2)} \cup P_j(\cb{u_3})$.
The first path of the second pair is $P_j(u_0)$.
Let $R_1 = Y_i(u_0) \cb{\cup Q_{m_i}(u_1)} \cup P_{m_i}(\cb{u_2})$,
$R_2 = P_j(\cb{u_2}) \cup Q_j(\cb{u_2}) \cup Q_j(\cb{u_3})$,
and $R_3 = X_i(u_0) \cb{\cup Q_{l_i}(u_1)} \cup X_{l_i}(\cb{u_2}) \cup X_{l_{l_i}}(\cb{u_3})$.
Then the second path of the second pair
is \cb{a path from $\xi_{v_1}(i)$ to $\xi_{v_2}(i)$ that is
a subgraph of} $R_1 \cup R_2 \cup R_3 \cup R_4 \cup R_5$, where
$R_4$ is a W6-path in the outer graph of \cb{$v_6$} joining the ends of $R_1$ and $R_2$,
and $R_5$ is a W6-path in the outer graph of \cb{$v_9$} joining the ends of $R_2$ and $R_3$.
See Figure~\ref{fig:abcd3}.
This completes main case 1.

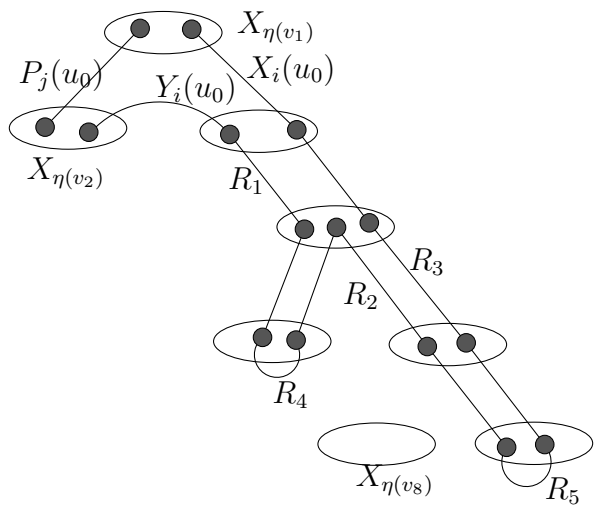
\begin{figure}[htb]
\centering
\input{image13.tikz}
\caption{Second pair when $l_i\in B-A$ and and $m_i\in A-B$.}
\label{fig:abcd3}
\end{figure}
\medskip

\noindent
{\bf Main case 2:} $A \cap B = \emptyset$.
It follows that $s$ is even and $|A| = |B| = s/2$.
Assume as a case that for some integer $i\in B$ either $l_i,m_i\in A$ or $l_i,m_i\in B$.
But the integers $l_i,m_i$ are pairwise distinct, and so if $l_i,m_i\in A$,
then there exists $j\in B$ such that $l_j,m_j\in B$,
and similarly if $l_i,m_i\in B$.
We may therefore assume that $l_i,m_i\in A$ and $l_j,m_j\in B$ for some distinct $i,j\in B$.
Let us recall that $u_2$ is the child of $v_5$
and $(v_5, v_6, v_7)$ is its trinity.
We let $\gamma$ map $t_0,t_1,t_2,t_3$ to $u_0,v_1,v_2,v_6$, respectively,
and we will prove that $t_0$ has property $AB_{ij}$ in $\eta'$.
To that end we need to construct two pairs of disjoint paths. The first pair is $Q_i(u_0)\cap Q_i(u_1)\cap P_i(u_2)$
and $Q_j(u_0)\cap Q_j(u_1)\cap P_j(u_2)$.
The first path of the second pair will consist of the union of
$X_i(u_0)$
with a subpath of $Q_{l_i}(u_1)$ from $X_{\eta(v_3)}$ to $X_{\eta(v_4)}$,
and $Y_i(u_0)$ with a subpath of $Q_{m_i}(u_1)$ from $X_{\eta(v_3)}$ to $X_{\eta(v_4)}$,
and a suitable W6-path in the outer graph of $v_4$
 joining their ends, and the second path will consist of the union of
$X_j(u_0)\cup \redsout{X_{l_j}(u_1)\cup X_{l_{l_j}}(u_2)}\cb{Q_{l_j}(u_1)\cup Q_{l_j}(u_2)}$ 
and $Y_j(u_0)\cup  \redsout{X_{\cb{m_j}}(u_1)\cup X_{m_{\cb{m_j}}}(u_2)}\cb{Q_{m_j}(u_1)\cup Q_{m_j}(u_2)}$
and a suitable W6-path in the outer a graph of $v_7$ joining their ends.
See Figure~\ref{fig:abcd4}.
This completes the case that for some integer $i\in B$ either $l_i,m_i\in A$ or $l_i,m_i\in B$.

\begin{figure}[htb]
\centering
\input{image14.tikz}
\caption{Second pair when $l_i,m_i\in A$ and $l_j,m_j\in B$ for some distinct $i,j\in B$.}
\label{fig:abcd4}
\end{figure}
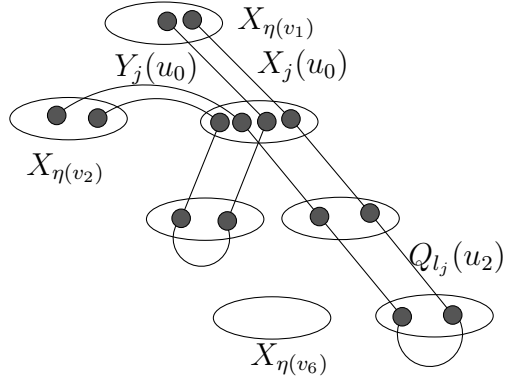

We may therefore assume that for every $i\in B$ one of $l_i,m_i$ belongs to $A$ and the other  belongs to $B$.
 Let us recall that for every $i\in B$ a subpath of $P_i(u_0)$ joins $\xi_{v_3}(l_i)$ to $\xi_{v_3}(m_i)$
in the outer graph  at $v_3$ and is disjoint from the $\eta$-torso at $u_0$, except for its ends.
Let $J$ be the union of these subpaths; then $J$ is a linkage from
$\{\xi_{v_3}(i)\,:\,i\in A\}$ to $\{\xi_{v_3}(i)\,:\,i\in B\}$.
For $i\in B$ the path $Q_i(u_0)$ is a subgraph of  the $\eta$-torso at $u_0$.
For $i\in A$ the intersection of the path $Q_i(u_0)$ with the $\eta$-torso at $u_0$  consists of two paths,
one from $X_{\eta(v_1)}$ to $X_{\eta(v_2)}$, and the other from $X_{\eta(v_2)}$ to $X_{\eta(v_3)}$.
Let $L$ denote the union of these subpaths over all $i\in A$.
It follows that $J\cup L\cup\bigcup_{i\in B}Q_i(u_0)$ is a linkage from $X_{\eta(v_1)}$ to $X_{\eta(v_2)}$,
and so by the minimality of the specified $u_0$-linkages,
it is equal to the specified left $u_0$-linkage.
It follows that $u_0$ has property $C$ in $\eta$.
\end{proof}

\begin{lemma}
\label{lemma4.16a}
Let $(T,X)$ be a tree-decomposition of a graph $G$ satisfying \cd{(W6) and (W7)}.
If there exists a regular cascade \cd{$\eta:T_3\emb T$} with orderings $\xi_t$
in which every major vertex has property \cd{$C$},
then there is a weak subcascade $\eta'$ of $\eta$ of height one such that
the major vertex in $\eta'$ has property \cd{$C_{ij}$} for some $i,j$.
\end{lemma}

\begin{proof}
Let the common confinement sets for $\eta$ be $A, B, C, D$.
For a major vertex $w\in V(T_3)$ with trinity $(v_1,v_2,v_3)$ there are disjoint paths in the $\eta$-torso at $w$
as in the definition of property C. For $a\in A$ and $b\in B$ let
 $R_a(w)$ denote the path with ends
$\xi_{v_1}(a)$ and $\xi_{v_2}(a)$, let
 $R_b(w)$ denote the path with ends
$\xi_{v_1}(b)$ and $\xi_{v_3}(b)$, and let
 $R_{ab}(w)$ denote the path with ends
$\xi_{v_2}(b)$ and $\xi_{v_3}(a)$.

Assume the major root of $T_3$ is $u_0$ and its trinity is $(v_1,v_2,v_3)$, and let $I$ be the
common intersection set of $\eta$.
Then $\eta(v_1),\eta(v_2),\eta(v_3)$ is a triad in $T$ with center $\eta(u_0)$ and for
all $i\in\{1,2,3\}$ we have $X_{\eta(v_i)}\cap X_{\eta(u_0)}=I=X_{\eta(v_1)}\cap X_{\eta(v_2)}\cap X_{\eta(v_3)}$,
and hence the triad is  not $X$-separable
%
by (W7). Thus by Lemma~\ref{cutbag} there is a path $R(u_0)$ connecting two of the three sets of disjoint paths in the $\eta$-torso at $u_0$. Assume
without loss of generality  that one end of
$R(u_0)$ is in a path  $R_i(u_0)$, where $i\in A$.
Then the other end of $R(u_0)$ is either in a path $R_j(u_0)$, where $j\in B$; or  in a path $R_{aj}(u_0)$,
where $j\in B$ and $a\in A$.
In the former case we define $a\in A$ to be such that $R_{aj}(u_0)$ is a path in the family.

%
Let the major root of $T_1$ be $t_0$ and its trinity be $(t_1, t_2, t_3)$. Let
$\gamma(t_0) = u_0$, $\gamma(t_1) = v_1$, $\gamma(t_2) = v_2$.
Let the major vertex that is the child of $v_3$ be $u_1$, and the trinity of
$u_1$ be $(v_3, v_4, v_5)$. Let $\gamma(t_3) = v_5$.
We will prove that $t_0$ has property $C_{ij}$ in $\eta' = \eta \circ \gamma$.
Let $b\in B$ be such that $R_{ib}(u_1)$ is a member of the family of the
disjoint paths in the $\eta$-torso at $u_1$
as in the definition of property C.
By Lemma~\ref{lem:W6}, there exists a $W6$-path $P$ in the
outer graph at $v_4$ joining $\xi_{v_4}(a)$ and $\xi_{v_4}(b)$.
By considering the paths $R_a(u_0)$, $R_j(u_0) \cup R_j(u_1)$,
$R_{aj}(u_0) \cup R_a(u_1) \cup P\cup R_{ib}(u_1)$ and $R(u_0)$
we find that $t_0$ has property $C_{ij}$ in $\eta'$, as desired.
\end{proof}

\begin{lemma}
\label{lemma4.17}
Let $s\geq 2$ be an integer and let $(T,X)$ be a tree-decomposition of a graph $G$ satisfying \cd{(W6)}.
Let \cd{$\eta:T_3\emb T$} be an ordered cascade in $(T,X)$ of height three and size $|I|+s$
with orderings $\xi_t$ and common intersection set $I$ such that every major vertex of $T_3$ has property $C_{ij}$
for some distinct $i, j\in\{1,2,\ldots,s\}$. Then there exists a weak subcascade
$\eta':T_{1}\emb T$  of $\eta$ of height one such that the unique major vertex of
 $T_1$ has property $B_{ij}$  in $\eta'$.
\end{lemma}

\begin{proof}
Assume that the three major vertices \cd{at height zero and one of $T_3$} are $u_0, u_1, u_2$.
Let the trinity at $u_0$ be $(v_1, v_2, v_3)$, the trinity at $u_1$ be $(v_2, v_4, v_5)$,
and the trinity at $u_2$ be $(v_3, v_6, v_7)$. Assume the major vertex of $T_1$ is $t_0$,
and its trinity is $(t_1, t_2, t_3)$. For a major vertex $w\in V(T_3)$ let  $R_i(w), R_j(w), R_{ij}(w)$
and $R(w)$ be as in the definition of property~$C_{ij}$.

We need to find a weakly  monotone \he\ $\gamma: T_{1}\emb T_3$ such that $\eta' = \eta \circ \gamma$
satisfies the requirement. Set $\gamma(t_0) = u_0$ and $\gamma(t_1) = v_1$.
Our choice for $\gamma(t_2)$ will be $v_4$ or $v_5$, depending on
which two of the three paths $R_i(u_1), R_j(u_1), R_{ij}(u_1)$ in the torso at $u_1$
the path $R(u_1)$ is connecting.
If $R(u_1)$ is between $R_i(u_1)$ and $R_j(u_1)$, then choose either $v_4$ or $v_5$
for $\gamma(t_2)$. If $R(u_1)$ is between $R_i(u_1)$ and $R_{ij}(u_1)$, then set
$\gamma(t_2) = v_4$, and if it is between $R_j(u_1)$ and $R_{ij}(u_1)$, then set
$\gamma(t_2) = v_5$. Do this similarly for $\gamma(t_3)$. Then $\eta' = \eta \circ \gamma$
will satisfy the requirement. In fact, we will prove this for the case when $R(u_1)$
is between $R_i(u_1)$ and $R_{ij}(u_1)$ and $R(u_2)$
is between $R_j(u_2)$ and $R_{ij}(u_2)$. See Figure~\ref{fig:dijtobij}. The other five cases are similar.

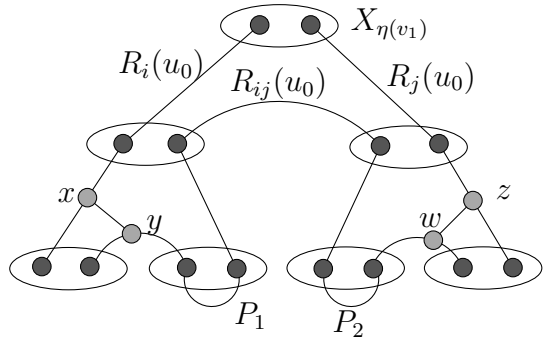
\begin{figure}[htb]
\centering
\input{image12.tikz}
\caption{The case when $R(u_1)$ is between $R_i(u_1)$ and $R_{ij}(u_1)$ and
$R(u_2)$ is between $R_j(u_2)$ and $R_{ij}(u_2)$.}
\label{fig:dijtobij}
\end{figure}

In this case, our choice is $\gamma(t_0) = u_0, \gamma(t_1) = v_1, \gamma(t_2) = v_4, \gamma(t_3) = v_7$.
Assume the two endpoints of $R(u_1)$ are $x$ and $y$ and the two endpoints of $R(u_2)$ are $w$ and $z$.
\cd{By Lemma \ref{lem:W6}, there exists a W6-path $P_1$
between $\xi_{v_5}(i)$ and $\xi_{v_5}(j)$ in the outer graph at $v_5$
and a W6-path $P_2$ between $\xi_{v_6}(i)$ and $\xi_{v_6}(j)$ in the outer graph at $v_6$.}
Now let $$P = yR_{ij}(u_1)\xi_{v_5}(i)\cup P_1 \cup R_j(u_1) \cup R_{ij}(u_0) \cup
R_i(u_2) \cup P_2 \cup \xi_{v_6}(j)R_{ij}(u_2)w,$$
$$L_i = R_i(u_0) \cup R_i(u_1) \cup R(u_1) \cup P \cup  wR_{ij}(u_2)\xi_{v_7}(i)$$ and
$$L_j = R_j(u_0) \cup R_j(u_2) \cup R(u_2) \cup P \cup yR_{ij}(u_1)\xi_{v_4}(j).$$
The tripods $L_i$ and $L_j$ show that
the major vertex of $\eta' = \eta \circ \gamma :T_{1}\emb T$ has property~$B_{ij}$.
\end{proof}

\begin{lemma}
\label{lemma4.15}
For every positive integers $h'$ and $w \geq 2$ there exists a positive integer $h = h(h', w)$
such that the following holds. Let $s$ be a positive integer such that $2 \leq s \leq w$.
Let $(T,X)$ be a tree-decomposition of a graph $G$ of width less than $w$ and satisfying \cd{(W6) and (W7)}.
Assume there exists a regular cascade $\eta:T_{h}\emb T$ of size $|I| + s$
with specified linkages that are minimal,
where $I$ is its common intersection set.
Then there exist distinct integers $i,j\in\{1,2,\ldots,s\}$ and
a weak subcascade $\eta':T_{h'}\emb T$  of $\eta$  of height $h'$ such that
\begin{itemize}
\item every major vertex of $T_{h'}$ has property $A_{ij}$ in $\eta'$, or
\item every major vertex of $T_{h'}$ has property $B_{ij}$ in $\eta'$
\end{itemize}
\end{lemma}

\begin{proof}
Let $h(a,k)$ be the function of Lemma~\ref{lemma4.10},
let $a_3=3h'$, $a_2=h(a_3,2{w\choose2})$, $a_1=\cb{5a_2}$ and $h=h(a_1,2)$.
Consider having property $C$ or not having property $C$ as colors, then by Lemma~\ref{lemma4.10}
there exists a monotone homeomorphic
embedding $\gamma: T_{a_1} \emb T_h$ such that
either $\gamma(t)$ has property $C$ in $\eta$ for every major vertex $t \in V(T_{a_1})$
or $\gamma(t)$ does not have property $C$ in $\eta$
for every major vertex $t \in V(T_{a_1})$.
By Lemma~\ref{lem:subcasc} $\eta_1 = \eta \circ \gamma: T_{a_1} \emb T$ is still a regular cascade
with specified linkages that are minimal.
Also, either $t$ has property $C$ in $\eta_1$ for every major vertex $t \in V(T_{a_1})$
or $t$ does not have property $C$ in $\eta_1$ for every major vertex $t \in V(T_{a_1})$.

If $t$ has property $C$ in $\eta_1$ for every major vertex $t\in V(T_{a_1})$, then
by Lemma \ref{lemma4.16a}
there exists a weak subcascade $\eta_2$ of $\eta_1$
of height $a_2$ such that every
major vertex of $T_{a_2}$ has property $C_{ij}$ in $\eta_2$ for some distinct $i,j \in \{1,2,...,s\}$.
Consider each choice of pair $i,j$ as a color; then
by Lemma \ref{lemma4.10}
there exists a  monotone homeomorphic embedding $\gamma_1: T_{a_3} \emb T_{a_2}$
such that for some distinct $i,j \in \{1,2,...,s\}$,
$\gamma_1(t)$ has property $C_{ij}$ in $\eta_2$ for every major vertex $t\in V(T_{a_3})$.
Let $\eta_3 = \eta_2 \circ \gamma_1$. Then by Lemma~\ref{lem:subcasc}
this  implies $t$ has property $C_{ij}$ in $\eta_3$ for every major vertex $t \in V(T_{a_3})$.
Then by Lemma \ref{lemma4.17}
there exists a weak subcascade $\eta_4: h' \emb a_3$ of $\eta_3$ such that every major vertex of $T_{h'}$
has property $B_{ij}$ in $\eta_4$. Hence $\eta_4$ is as desired.

If $t$ does not have property $C$ in $\eta_1$ for every major vertex $t\in V(T_{a_1})$,
then by Lemma~\ref{lemma4.14}
there exists a weak subcascade $\eta_2$ of $\eta_1$ of height $a_2$ such that
every major vertex of $T_{a_2}$ has property $A_{ij}$ or $B_{ij}$ for some distinct
$i,j \in \{1,2,...,s\}$. Consider each property $A_{ij}$ or $B_{ij}$ as a color; then
by Lemma \ref{lemma4.10}
there exists a monotone homeomorphic embedding $\gamma_1: T_{h'} \emb T_{a_2}$
such that for some distinct $i,j\in \{1,2,...,s\}$, either
$\gamma_1(t)$ has property $A_{ij}$ in $\eta_2$ for every major vertex $t\in V(T_{h'})$
or $\gamma_1(t)$ has property $B_{ij}$ in $\eta_2$ for every major vertex $t\in V(T_{h'})$.
Let $\eta_3 = \eta_2 \circ \gamma_1$.
Then $t$ has property $A_{ij}$ in $\eta_3$ for every major vertex $t \in V(T_{h'})$
or $t$ has property $B_{ij}$ in $\eta_3$ for every major vertex $t \in V(T_{h'})$ by Lemma~\ref{lem:subcasc}.
Hence $\eta_3$ is as desired.
\end{proof}

\section{Proof of Theorem~\ref{thm3}}
\label{sec4}

By Lemmas~\ref{lemma2.1} and~\ref{lemma2.2}
Theorem~\ref{thm3} is equivalent to the following theorem.
\begin{theorem}
\label{thm4}
For any positive integer $k$, there exists a positive integer $p = p(k)$
such that for every $2$-connected graph $G$,
if $G$ has path-width at least $p$, then $G$ has  a minor isomorphic to $\ca{\mathcal{P}_k}$ or $\ca{\mathcal{Q}_k}$.
\end{theorem}

We need the following three lemmas. 
We state them in greater generality than immediately necessary in order to be able to use them 
in a subsequent paper.

\begin{definition}
For a positive integer $k$,
let $\mathcal{P}'_k$ be the graph consisting of 	$CT_k$
plus two distinct vertices each adjacent to every leaf of $CT_k$,
let $\mathcal{P}''_k$ be the graph consisting of 	$CT_k$
plus three distinct vertices each adjacent to every leaf of $CT_k$,
let $\mathcal{Q}'_k$ be $\mathcal{Q}_k$ plus a vertex, called the {\em apex}, adjacent to its vertices of degree two, and
let $\mathcal{Q}''_k$ be $\mathcal{Q}_k$ plus two vertices, called the {\em apices}, each adjacent to its vertices of degree two.
\end{definition}

\begin{lemma}
\label{lem:A2P}
Let $(T,X)$ be a \td\ of a graph $G$ satisfying {\rm (W6)}, let $k\ge1$ be an integer, let $h=2k+1$,
let $\eta:T_h\emb T$ be an ordered cascade in $(T,X)$ with orderings $\xi_t$
of height $h$ and size $s+|I|$, where $I$ is the
common intersection set, and let $i,j\in\{1,2,\ldots,s\}$ be distinct and
 such that every major vertex of $T_h$ has property $A_{ij}$ in~$\eta$.
Then
\begin{itemize}
\item[\rm(i)] $G$ has a minor isomorphic to $\mathcal{P}_h$,
\item[\rm(ii)] if $|I|\ge1$, then $G$ has a minor isomorphic to $\mathcal{P}_h'$, and
\item[\rm(iii)] if $|I|\ge2$, then $G$ has a minor isomorphic to $\mathcal{P}_h''$.
\end{itemize}
\end{lemma}

\begin{proof}
We only prove (iii), since the other two statements are analogous, and, in fact, easier.
Let $R$ be the union of the corresponding tripods as in the definition of property $A_{ij}$, over all
major vertices $t\in V(T_{h})$ at height at most $h-2$. It follows that $R$ is the union
of two disjoint trees, each containing a subtree isomorphic to $CT_{k}$.
Let $t$ be a minor vertex of $T_{h}$ at height $h-1$, and let $x,y\in I$.
By Lemma~\ref{lem:W6} there exist W6-paths with ends $u$ and $v$
in the outer graph at $t$ for every pair of distinct vertices $u,v\in\{\xi_t(i),\xi_t(j),x,y\}$.
By contracting one of the trees comprising   $R$ and by considering these W6-paths we deduce that $G$
has a $\mathcal{P}''_{k}$ minor, as desired.
\end{proof}

\begin{lemma}
\label{lem:B2Q}
Let $(T,X)$ be a \td\ of a graph $G$ satisfying {\rm (W6)}, let $h\ge1$ be an integer,
let $\eta:T_h\emb T$ be an ordered cascade in $(T,X)$ with orderings $\xi_t$
of height $h$ and size $s+|I|$, where $I$ is the
common intersection set, and let $i,j\in\{1,2,\ldots,s\}$ be distinct and
 such that every major vertex of $T_h$ has property $B_{ij}$ in $\eta$.
Let $t$ be the minor root of $T_h$, and
let $w_1w_2$ be the base edge of $\mathcal{Q}_h$. Then 
\begin{itemize}
\item[\rm(i)]
$G$ has a minor isomorphic
to $\mathcal{Q}_h-w_1w_2$ in such a way that $\xi_t(i)$ belongs to the node of $w_1$ and
 $\xi_t(j)$ belongs to the node of $w_2$,
\item[\rm(ii)]if $x\in I$, then
$G$ has a minor isomorphic
to $\mathcal{Q}'_{h-1}-w_1w_2$ in such a way that $\xi_t(i)$ belongs to the node of $w_1$,
 $\xi_t(j)$ belongs to the node of $w_2$ and $x$ belongs to the node of the apex, and
\item[\rm(iii)]if $x,y\in I$ are distinct, then
$G$ has a minor isomorphic
to $\mathcal{Q}''_{h-1}-w_1w_2$ in such a way that $\xi_t(i)$ belongs to the node of $w_1$,
 $\xi_t(j)$ belongs to the node of $w_2$ and $x$ and $y$ belong to the nodes of the apices, respectively.
\end{itemize}
\end{lemma}

\begin{proof}
We only prove (i), noting that the other two statements follow similarly as the previous lemma.
We proceed by induction on $h$.
Let $t_0$ be the major root of $T_h$,  let $(t_1,t_2,t_3)$ be its trinity, and
let $L_i$ and $L_j$ be the tripods in the $\eta$-torso at $t_0$ as in the definition of property $B_{ij}$.
The graph $L_i\cup L_j$  contains a path $P$ joining $\xi_{t_1}(i)$ to $\xi_{t_1}(j)$, which
shows that the lemma holds for $h=1$.

We may therefore assume that $h>1$ and that the lemma holds for $h-1$.
For $k\in\{2,3\}$
let $R_k$ be the subtree of $T_{h}$ rooted at $t_k$,
let $\eta_k$ be the restriction of $\eta$ to $R_k$, and let $G_k$ be the subgraph of $G$ induced by $\bigcup\{X_r:r\in sp(\eta_k)\}$.
By the induction hypothesis applied to $\eta_k$ and $G_k$, the graph $G_k$ has a minor isomorphic to
$\mathcal{Q}_{h-1}-u_1u_2$ in such a way that $\xi_{t_k}(i)$ belongs to the node of $u_1$ and
 $\xi_{t_k}(j)$ belongs to the node of $u_2$, where $u_1u_2$ is the base edge of $\mathcal{Q}_{h-1}$.
By using these  two minors, the path $P$ and the rest of the tripods $L_i$ and $L_j$
we find that $G$ has the desired minor.
\end{proof}

\begin{lemma}
\label{lem:AB2PQ}
For every two integers $w,k\ge1$ there exists an integer $h\ge2k+1$ such that the following holds.
Let $(T,X)$ be a \td\ of a graph $G$ of width less than $w$ satisfying {\rm (W6)} and {\rm(W7)}, and
let $\eta:T_h\emb T$ be a regular cascade in $(T,X)$ 
of height $h$ and size $s+|I|$ with specified $t_0$-linkages that are
minimal for every major vertex $t_0\in V(T_h)$, where $s\ge2$ and $I$ is the
common intersection set.
Then
\begin{itemize}
\item[\rm(i)] $G$ has a minor isomorphic to $\mathcal{P}_k$ or $\mathcal{Q}_k$,
\item[\rm(ii)] if $|I|\ge1$, then $G$ has a minor isomorphic to $\mathcal{P}_k'$ or $\mathcal{Q}_k'$, and
\item[\rm(iii)] if $|I|\ge2$, then $G$ has a minor isomorphic to $\mathcal{P}_k''$ or $\mathcal{Q}_k''$.
\end{itemize}
\end{lemma}

\begin{proof}
Let $h'=2k+1$ and let $h = h(h', w)$ be the number as in Lemma~\ref{lemma4.15}.
We claim that $h$ satisfies the conclusion of the lemma.
To prove that let $G,(T,X)$ and $\eta$ be as in the statement of the lemma.
By Lemma~\ref{lemma4.15}
there exist distinct integers $i,j\in\{1,2,\ldots,s\}$ and
a weak subcascade  $\eta':T_{h'}\emb T$  of $\eta$  of height $h'$ such that
\begin{itemize}
\item every major vertex of $T_{h'}$ has property $A_{ij}$ in $\eta'$, or
\item every major vertex of $T_{h'}$ has property $B_{ij}$ in $\eta'$
\end{itemize}
If every major vertex of $T_{h'}$ has property $A_{ij}$ in $\eta'$, then
by Lemma~\ref{lem:A2P}  $G$
has a minor isomorphic to $\mathcal{P}_{k}$, as desired.
We may therefore assume that  every major vertex of $T_{h'}$ has property $B_{ij}$ in $\eta'$.
It follows from Lemma~\ref{lem:B2Q} that $G$
has a minor isomorphic to $\mathcal{Q}_{h'-1}$, as desired.
This proves (i).
The other two statements follow similarly.
\end{proof}

We deduce Theorem~\ref{thm4} from the following lemma.

\begin{lemma}
\label{lemma4.18}
Let $k$ and $w$ be positive integers. There exists a number $p= p(k,w)$
such that for every $2$-connected graph $G$, if $G$ has tree-width less than $w$
and path-width at least $p$, then $G$ has a minor isomorphic to $\ca{\mathcal{P}_k}$ or $\ca{\mathcal{Q}_k}$.
\end{lemma}

\begin{proof}
Let $h = h(k, w)$ be the number as in Lemma~\ref{lem:AB2PQ},
and let $p$ be as in Theorem~\ref{thm:regular} applied to $a=h$ and $w$. 
We claim that $p$ satisfies the conclusion of the lemma.
By Theorem~\ref{thm:regular}, there exists a tree-decomposition $(T,X)$ of $G$ such that:
\begin{itemize}
\item $(T, X)$ has width less than $w$,
\item $(T, X)$ satisfies (W1)--(W7), and
\item for some $s$, where $2 \leq s \leq w$, there exists a regular cascade $\eta:T_h\emb T$
of height $h$ and size $s$ in $(T, X)$ with specified $t_0$-linkages that are
minimal for every major vertex $t_0\in V(T_h)$.
\end{itemize}
Let $I$ be the common intersection set of $\eta$, let $\xi_t$ be the orderings, and let $s_1 = s - |I|$. Then $s_1\ge1$ by
the definition  of injective cascade.

Assume first that $s_1 = 1$. Since $s\ge2$, it follows that $I\ne\emptyset$. Let $x\in I$.
Let  $R$ be  the union of the left and right specified $t$-linkage with respect to $\eta$, over all
major vertices $t\in V(T_h)$ at height at most $h-2$.
The minimality of the specified linkages implies that $R$ has a subtree
 isomorphic to a subdivision of $CT_{\lfloor(h-1)/2\rfloor}$.
Let $t$ be a minor vertex of $T_h$ at height $h-1$.
By Lemma~\ref{lem:W6} there exists a W6-path with ends $\xi_t(1)$ and $x$ 
and every internal vertex in the outer graph at $t$.
Since $h\ge2k+1$,
the union of  $R$ and these W6-paths shows that $G$ has a $\mathcal{P}_{k}$ minor, as desired.

We may therefore assume that  $s_1 \geq 2$. 
The lemma now follows from  Lemma~\ref{lem:AB2PQ}.
%
\end{proof}


\begin{proof}[Proof of Theorem~\ref{thm4}]
Let a positive integer $k$ be given.
By Theorem~\ref{thm1} there exists an integer $w$ such that every graph of
tree-width at least $w$ has a minor isomorphic to~$\ca{\mathcal{P}_k}$.
Let $p=p(k,w)$ be as in Lemma~\ref{lemma4.18}.
We claim that $p$ satisfies the conclusion of the theorem.
Indeed,
let $G$ be a $2$-connected graph of path-width at least $p$.
 By Theorem~\ref{thm1}, if $G$ has tree-width at least $w$,
then $G$ has a minor isomorphic to $\ca{\mathcal{P}_k}$, as desired.
We may therefore
assume that the tree-width of $G$ is less than $w$.
By Lemma~\ref{lemma4.18}
 $G$ has a minor isomorphic to $\ca{\mathcal{P}_k}$ or $\ca{\mathcal{Q}_k}$, as desired.
\end{proof}

\section*{Acknowledgments}
The second author would like to acknowledge that Christopher Anhalt wrote a
``Diplomaarbeit" under his supervision in 1994 on the subject of this paper,
but the proof is incomplete.
The results of this paper form part of the Ph.D.\ dissertation~\cite{DanPhD}
of the first author, written under the guidance of the second author.
After our paper has been posted on the arXiv and submitted, a
different proof by Huynh, Joret, Micek and Wood~\cite{huyjormicwoo} appeared.
Our proof has the advantage that it can be generalized to graphs of higher connectivity~\cite{DanPhD}.


\end{document}

%% file: w7insert.tex
\section{Linked Tree-decompositions}
\label{sec:linked}

In this section we review properties of \td s established in
\cite{OpoOxlTho, thomas90}, and state our main lemma.
The proof of the following easy lemma can be found, for instance,
in
\cite{thomas90}.

\begin{lemma}
\label{cutbag}
Let $(T,Y)$ be a tree-decomposition of a graph $G$, and let $H$
be a
connected subgraph of $G$ such that $V(H)\cap Y_{t_1} \ne\emptyset\ne
V(H)\cap
Y_{t_2}$, where $t_1,t_2\in V(T)$.  Then $V(H)\cap Y_t\ne\emptyset$ for
every
$t\in V(T)$ on the path between $t_1$ and $t_2$ in $T$.
\end{lemma}

\noindent A tree-decomposition $(T,Y)$ of a graph $G$ is said to be
{\em linked} if

\begin{itemize}
\item[(W3)]for every two vertices $t_1,t_2$ of $T$ and every positive integer
$k$,
either there are $k$ disjoint paths in $G$ between $Y_{t_1}$ and
$Y_{t_2}$, or there is a vertex $t$ of $T$ on the path between
$t_1$ and $t_2$ such that $|Y_t| < k$.
\end{itemize}

\noindent It is worth noting that, by Lemma~\ref{cutbag}, the two alternatives in (W3)
are
mutually exclusive.  The following is proved in \cite{thomas90}.

\begin{lemma}
\label{linked}
If a graph $G$ admits a tree-decomposition of width at most $w$,
where $w$ is some integer, then $G$ admits a linked tree-decomposition
of
width at most $w$.
\end{lemma}

  Let $(T,Y)$ be a \td\ of a graph $G$,
let $t_0\in V(T)$, and let $B$ be a component of $T\backslash
t_0$.
We say that a vertex $v\in Y_{t_0}$ is $B$-{\em tied} if $v\in Y_t$
for some $t\in V(B)$.  We say that a path $P$
in $G$ is $B$-{\em confined} if $|V(P)|\ge3$ and every internal vertex
of $P$ belongs
to  $\bigcup\limits_{t\in V(B)} Y_t - Y_{t_0}$.  We wish to consider the
following three properties of $(T,Y)$:

\begin{itemize}
\item[(W4)]if $t,t'$ are distinct vertices of $T$, then $Y_t\ne Y_{t'}$,

\item[(W5)]if $t_0\in V(T)$ and $B$ is a component of
$T\backslash t_0$, then
$\bigcup\limits_{t\in V(B)} Y_t - Y_{t_0} \ne\emptyset$,

\item[(W6)]if $t_0\in V(T)$, $B$ is a component of
$T\backslash t_0$, and $u,v$ are $B$-tied
vertices in $Y_{t_0}$, then  there is a $B$-confined path in $G$
between $u$ and $v$.
\end{itemize}

\noindent The following strengthening of Lemma~\ref{linked} is proved
in~\cite{OpoOxlTho}.

\begin{lemma}
\label{lem:w1w6}
If a graph $G$ has a tree-decomposition of width at most
$w$, where $w$ is some integer,
then it has a tree-decomposition of width at most $w$ satisfying {\rm(W1)--(W6)}.
\end{lemma}

We need one more condition, which we now introduce.
Let $T$ be a tree. If $t_1,t_2\in V(T)$, then by $T[t_1,t_2]$ we denote the
vertex-set of the unique path in $T$  with ends $t_1$ and $t_2$. 
A {\em triad in $T$} is a triple $t_1,t_2,t_3$ of
vertices of $T$ such that there exists a vertex $t$ of $T$, called the
{\em center}, such that $t_1,t_2,t_3$ belong to different components
of $T\backslash t$.
Let $(T,W)$ be a tree-decomposition of a graph $G$, and let
$t_1,t_2,t_3$ be a triad in $T$ with center $t_0$. 
The {\em torso of $(T,W)$ at $t_1,t_2,t_3$}
is the subgraph of $G$ induced by the set $\bigcup W_t$,
the union taken over all vertices $t\in V(T)$ such that
either $t\in\{t_1,t_2,t_3\}$, or for all
$i\in\{1,2,3\}$, the vertex $t$ belongs to the component of $T\backslash t_i$ containing $t_0$.
We say that the triad $t_1,t_2,t_3$
is {\em$W$-separable} if, letting $X=W_{t_1}\cap W_{t_2}\cap W_{t_3}$,
the graph obtained from
the torso of $(T,W)$ at $t_1,t_2,t_3$ by deleting $X$
 can be partitioned into three disjoint non-null graphs
$H_1,H_2,H_3$ in such a way that for all distinct $i,j\in\{1,2,3\}$,
all $k\in \{1,2,3\}$ and all $t\in T[t_j,t_0]$,
$|V(H_i)\cap W_t|\ge|V(H_i)\cap W_{t_j}|=|W_{t_k}-X|/2\ge 1$.
(Let us remark that this condition implies that
$|W_{t_1}|=|W_{t_2}|=|W_{t_3}|$  and $V(H_i)\cap W_{t_i}=\emptyset$ for
$i=1,2,3$.)
The last property of a \td\ $(T,W)$ that we wish to consider is

\begin{itemize}
\item[(W7)]if $t_1,t_2,t_3$ is a $W$--separable triad in $T$ with center $t$, then
there exists an integer $i\in\{1,2,3\}$ with
$W_{t_i}\cap W_{t}-\left(W_{t_1}\cap W_{t_2}\cap W_{t_3}\right)
\not=\emptyset$.
\end{itemize}

\noindent The following is our main lemma.

\begin{theorem}
\label{w7}
If a graph $G$ has a tree-decomposition of width at most
$w$, where $w$ is some integer,
then it has a tree-decomposition of width at most $w$ satisfying
{\rm(W1)--(W7)}.
\end{theorem}

\section{A Quasi-order on Trees}
\label{sec:quasi}

A \dfn{quasi-ordered set} is a pair $(Q,\le)$, where $Q$ is a set
and $\le$ is a \dfn{quasi-order}; that is, a reflexive and transitive
relation on $Q$.  If $q,q'\in Q$ we define $q<q'$ to mean that
$q\le q'$ and $q'\not\le q$.  We say that $q,q'$ are \dfn{$\le$-equivalent}
if $q\le q'\le q$. We say that $(Q,\le)$ is a {\em linear quasi-order}
if for every two elements $q,q'\in Q$ either $q\le q'$ or $q'\le q$ or both.
Let $(Q,\le)$ be a linear quasi-order.
 If $A,B\subseteq Q$ we say that
$B\le$-\dfn{dominates} $A$ if
the elements of $A$ can
be listed as $a_1\ge a_2\ge\cdots\ge a_k$ and the elements of
$B$ can be listed as $b_1\ge b_2\ge \cdots\ge b_l$, and there exists
an integer $p$ with $1\le p\le  \min\{k,l\}$ 
such that $a_i\le b_i\le a_i$
for all $i=1,2,\dots, p$, and either 
\cb{$p< \min\{k,l\}$ and $a_{p+1}<b_{p+1}$, or $p= k$ and $k\leq l$}.

\begin{lemma}
\label{lem:linqo}
If $(Q,\le)$ is a linear quasi-order,
then $\le$-domination is a linear quasi-order on  the set of subsets of $Q$.
\end{lemma}

\begin{proof}
It is obvious that $\le$-domination is reflexive.
Assume that $B$ $\le$-dominates $A$ and $C$ $\le$-dominates $B$.
Assume that the elements of $A$ can be listed as $a_1\ge a_2\ge\cdots\ge a_k$,
the elements of $B$ can be listed as $b_1\ge b_2\ge \cdots\ge b_l$, and
the elements of $C$ can be listed as $c_1\ge c_2\ge \cdots\ge c_m$.
By definition, there exists an integer $p_1$ with $1 \leq p_1 \leq \min\{k,l\}$
such that $a_i\le b_i\le a_i$
for all $i=1,2,\dots, p_1$, and either $p_1< \min\{k,l\}$ and $a_{p_1+1}<b_{p_1+1}$, or
$p_1 = k \leq l$;
and there exists an integer $p_2$ with $1 \leq p_2\leq \min\{l,m\}$
such that $b_i\le c_i\le b_i$
for all $i=1,2,\dots, p_2$, and either $p_2 < \min\{l,m\}$ and $b_{p_2+1}<c_{p_2+1}$, or
$p_2 = l \leq m$. Let $p = \min\{p_1, p_2\}$.
Then $a_i \leq c_i \leq a_i$ for all $i=1,2,\ldots,p$.
If either $p_1< \min\{k,l\}$ and $a_{p_1+1}<b_{p_1+1}$, or
$p_2 < \min\{l,m\}$ and $b_{p_2+1}<c_{p_2+1}$, then $p < \min\{k,m\}$ and $a_{p+1} < c_{p+1}$.
If $p_1 = k \leq l$ and $p_2 = l \leq m$, then $p = k \leq m$.
Therefore, $C$ $\le$-dominates $A$, and so $\le$-domination is transitive.

Now let $A, B$ be as above, and let $p$ be the maximum integer such that $p\le \min\{k,l\}$ and
$a_i\le b_i\le a_i$ for all $i=1,2,\ldots,p$. Then if $p<\min\{k,l\}$, then
$A$ $\le$-dominates $B$ if $a_{p+1} > b_{p+1}$ and
$B$ $\le$-dominates $A$ if $a_{p+1} < b_{p+1}$. If $p = \min\{k,l\}$
then $A$ $\le$-dominates $B$ if $k \geq l$ and $B$ $\le$-dominates $A$ if $k \leq l$.
Hence, $\le$-domination is linear.
\end{proof}

We say that $B$ \dfn{strictly $\le$-dominates} $A$ if $B$
$\le$-dominates $A$ in such a way that the numberings and integer $p$
can be chosen in such a way that either 
$p< \min\{k,l\}$, or 
$p=k$ and $k<l$.

\begin{lemma}
\label{lem:domin}
Let $(Q,\le)$ be a linear quasi-order,  let  $A,B\subseteq Q$, and let
$B$ $\le$-dominate $A$. Then  $B$ strictly $\le$-dominates $A$
if and only if  $A$ does not $\le$-dominate $B$.
\end{lemma}

\begin{proof}
Let $p$ be as in the definition of $B$ $\le$-dominates $A$.
Then $p< \min\{k,l\}$ and $a_{p+1}<b_{p+1}$, or $p= k \leq l$.
Assume $B$ strictly $\le$-dominates $A$.
If $p < \min\{k,l\}$
then $a_{p+1} < b_{p+1}$, so $A$ does not $\le$-dominate $B$. If $p=k<l$ then
$A$ also does not $\le$-dominate $B$. Conversely, if $A$ does not $\le$-dominate $B$,
then $p < \min\{k,l\}$ or $k<l$, so $B$ strictly $\le$-dominates $A$.
\end{proof}

Let $G$ be a graph and let $P$ be a subgraph of $G$. By a $P$-bridge of $G$ we mean a
subgraph $J$ of $G$ such that either
\begin{itemize}
\item $J$ is isomorphic to the complete graph on two vertices with $V(J)\subseteq V(P)$
and $E(J)\cap E(P)=\emptyset$, or
\item $J$ consists of a component of $G-V(P)$ together with all edges from that component to $P$.
\end{itemize} 

We now define a linear  quasi-order $\le$ on the class of finite trees as
follows. Let $n\ge 1$ be an integer, and suppose that
$T\le T'$ has been defined for all trees $T$ on fewer than $n$
vertices.  Let $T$ be a tree on $n$ vertices, and let $T'$
be an arbitrary tree. We define $T\le T'$ if either $|V(T)|<|V(T')|$,
or $|V(T)|=|V(T')|$ and for every maximal path $P'$ of $T'$ there exists a
maximal path
$P$ of $T$ such that the set of $P'$-bridges of $T'$ $\le$-dominates
the set of $P$-bridges of $T$.
It follows from Lemma~\ref{spine} below that $\le$ is indeed a linear quasi-order;
in particular, it is well-defined.


If $T,T'$ are trees, $P$ is a path in $T$ and $P'$ is a path in $T'$
we define $(T,P)\preceq (T',P')$ if either $|V(T)|<|V(T')|$, or
$|V(T)|=|V(T')|$ and the set of $P'$-bridges of $T'$ $\le$-dominates
the set of $P$-bridges of $T$.

\begin{lemma}
\label{spine}
(i) For every tree $T$ there exists a maximal path $P(T)$ in $T$ such that
$(T,P(T))\preceq(T,P)$ for every maximal path $P$ in $T$.\\
(ii) For  every two trees $T,T'$, we have $T\le T'$ if and only if $(T,P(T))\preceq
(T',P(T'))$.\\
(iii) The ordering $\le$ is a linear quasi-order on the class of finite trees.
\end{lemma}

\begin{proof}
We prove all three statements simultaneously by induction. Let $n\ge1$ be an
integer,  assume inductively that all three statements have been proven for
trees on fewer than $n$ vertices, and let $T$ be a tree on $n$ vertices.\\
(i) Statement (i) clearly holds for one-vertex trees, and so we may assume that
$n\ge2$.
Let $\cal B$ be the set of all $P$-bridges of $T$ for all maximal paths $P$ of $T$.
Then every member of $\cal B$ has fewer than $n$ vertices,  and hence $\cal B$
is a linear quasi-order by $\le$ by the induction hypothesis applied to (iii).
By Lemma~\ref{lem:linqo} the set of subsets of $\cal B$ is linearly quasi-ordered
by $\le$-domination.
It follows that there exists a maximal path $P(T)$ in $T$ such that the set
of $P(T)$-bridges of $T$ is minimal under $\le$-domination.\\
(ii) The statement is obvious when $|V(T)| \neq |V(T')|$, so assume $n=|V(T)| = |V(T')|$,
and let $\cal B$ be the set of all $P$-bridges of $T$ for all maximal paths $P$ of $T$
and  the set of all $P'$-bridges of $T'$ for all maximal paths $P'$ of $T'$.
Then as in (i) the subsets of $\cal B$ are linearly quasi-ordered
by $\le$-domination.
If $T \leq T'$, then by definition there exists a maximal path $P$ of $T$ such that
$(T,P) \preceq (T', P(T'))$. Hence $(T,P(T))) \preceq (T', P(T'))$ follows from (i).
If $(T,P(T))) \preceq (T', P(T'))$, then by (i) $(T,P(T)))\preceq (T', P')$
for every maximal path $P'$ in $T'$, so $T \leq T'$.\\
(iii) Let $T$ and $T'$ be two trees. We may assume that $n=|V(T)| = |V(T')|$.
Let $\cal B$ be as in (ii); then subsets of $\cal B$ are linearly quasi-ordered
by $\le$-domination. Then
 either $(T,P(T))\preceq (T',P(T'))$ or
$(T',P(T'))\preceq (T,P(T))$,  and so by (ii) $\le$ is linear.
\end{proof}

For a tree $T$, the path $P(T)$ from Lemma~\ref{spine}(i) will be called a
\dfn{spine of} $T$.
For later application we need the following lemma.

\begin{lemma}
\label{ranklemma} Let $T,T'$ be trees on the same number of vertices,
let $P'$ be a spine of $T'$, and let $P$ be a path in $T$.
If the set of $P'$-bridges of $T'$ strictly $\le$-dominates
the set of $P$-bridges of $T$, then $T<T'$.
\end{lemma}

\begin{proof}
We have
 $(T,P)\preceq (T',P')$ and $(T',P')\not\preceq (T,P)$ by Lemma~\ref{lem:domin}.
Let $P_1$ be a maximal path that contains $P$; then $(T,P_1)\preceq (T, P)$.
Therefore, $(T,P_1) \preceq (T', P')$ and $(T',P')\not\preceq (T,P_1)$. 
By Lemma~\ref{spine}(i), $(T,P(T)) \preceq (T,P_1)\preceq(T', P')$ and $(T',P')\not\preceq (T,P(T))$.
By Lemma \ref{spine}(ii), $T \leq T'$ and $T' \not \leq T$. Therefore, $T < T'$.
\end{proof}

By a \dfn{rank} we mean a class of $\le$-equivalent trees.
If $r$ is a rank we say that \dfn{$T$ has rank $r$} or that the \dfn{rank of
$T$ is $r$} if $T\in r$. The class of all ranks will be denoted by
$\sR$.

Let $T$ be a tree, and let $t$ be a vertex of $T$.
By a spine-decomposition of $T$ relative to $t$ we mean a
sequence $(T_0,P_0, T_1,P_1,\dots, T_l, P_l)$ such that
\begin{itemize}
\item[(i)] $T_0=T$,
\item[(ii)] for $i=0,1,\dots, l$, $P_i$ is a spine of $T_i$, and
\item[(iii)] for $i=1,2,\dots, l$, $t\notin V(P_{i-1})$ and $T_i$ is
the $P_{i-1}$-bridge of $T_{i-1}$ containing $t$.
\end{itemize}

\begin{lemma}
\label{spine2}Let $T$ be a tree, let $t$ be a vertex of $T$ of degree three
with neighbors \cb{$t'_1,t'_2,t'_3$}, and let $(T_0, P_0, T_1, P_1,\dots,
T_l, P_l)$ be a spine-decomposition of $T$ relative to $t$
with $t\in V(P_l)$.
Then exactly two of \cb{$t'_1,t'_2,t'_3$} belong to $V(P_l)$, say $t'_1$ and $t'_2$.
Let $r_3,r'_3$ be adjacent vertices of $T$ such that
\cb{$r_3,r'_3,t'_3,t$} occur on a path of $T$ in the order listed.
Thus possibly $t_3'=r'_3$, but  $t_3'\ne r_3$.
Let $T'$ be obtained from $T$ by subdividing the edge $r_3r_3'$ twice (let $r''_3,r'''_3$
be the new vertices so that \cb{$r_3',r_3'',r_3''',r_3$} occur on a path of $T'$ in the order listed),
deleting the edge $tt_1'$, contracting the edges $tt_2'$ and $tt_3'$
and adding an edge joining $t_1'$ and $r'''_3$.
Then $T'$ has strictly smaller rank than $T$.
\end{lemma}

\begin{proof}
Let $T'_0=T'$ and for $i=1,2,\dots, l$, let $T'_i$ be the
$P_{i-1}$-bridge of $T_{i-1}'$ containing $r_3'''$.
Let $P'$ be the unique maximal path in $T'$ with $V(P_l)-\{t,t_2'\}\cup\{r_3'\}\subseteq V(P')$.
From the definition of a spine-decomposition and the fact that $t_3'\not\in V(P_l)$ we deduce
that $r_3\in V(T_i)$ for all $i=0,1,\ldots,l$.
It follows that $r_3\in V(T_i')$ and
$|V(T_i)|  = |V(T'_i)|$ for all $i=0,1,\ldots,l$. The
$P_l$-bridge of $T_l$ that contains $r_3$ is replaced by $P'$-bridges of $T'_l$
with smaller cardinalities. Other $P_l$-bridges of $T_l$ are unchanged in $T'$.
Therefore, the set of $P_l$-bridges of $T_l$
strictly $\le$-dominates the set of $P'$-bridges of $T'_l$,
and hence $T'_l<T_l$ by Lemma~\ref{ranklemma}.
This implies, by induction on $l-i$ using Lemma~\ref{ranklemma}, that $T'_i<T_i$ for all $i=0,1,\ldots,l$; that is,
$T'$ has smaller rank than $T$.
\end{proof}

\section{A Theorem about Tree-decompositions}
\label{sec:thm}


Let $(T,Y)$ be a tree-decomposition of a graph $G$,
let $n$ be an integer,
and let $r$ be a rank.  By an \dfn{$(n,r)$--cell}
in $(T,Y)$ we mean any component of the restriction of $T$ to $\{t\in
V(T):|Y_t|\ge n\}$ that has rank at least $r$.  Let us remark that
if $K$
is an $(n,r)$-cell in $(T,Y)$ and $r\ge r'$, then $K$ is an
$(n,r')$-cell as well.
The {\it size} of a tree-decomposition $(T,Y)$ is the family of numbers
$$
(a_{n,r}:n \ge 0, r\in\sR),\leqno(1)
$$
where $a_{n,r}$ is the number of $(n,r)$-cells in $(T,Y)$.  Sizes are
ordered
lexicographically; that is, if
$$
(b_{n,r}: n\ge 0, r\in\sR)\leqno(2)
$$
is the size of another tree-decomposition $(R,Z)$ of the graph $G$, we
say that
(2) is {\em smaller than} (1) if there are an integer $n\ge0$ and a rank
 $r\in\sR$ such that $a_{n,r}
> b_{n,r}$ and $a_{n',r'} = b_{n',r'}$ whenever either $n'>n$, or $n'=n$ and
 $r' > r$.

\begin{lemma}
\label{wellorder}The relation ``to be smaller than" is a well--ordering on the
set of sizes
of tree--decompositions of $G$.
\end{lemma}

\begin{proof} 
Since this ordering is clearly linear, it is enough to show that it is
well--founded.
Suppose for a contradiction that
${\{ (a^{(i)}_{n,r}: n\ge 0, r\in\sR)\}}_{i=1}^\infty$
is a strictly decreasing sequence of sizes, and for $i=1,2,\ldots, $ let $n_i,
r_i$ be such that
$a^{(i)}_{n_i, r_i} > a^{(i+1)}_{n_i,r_i}$ and
$a^{(i)}_{n, r} = a^{(i+1)}_{n,r}$ for $(n,r)$ such that either $n>n_i$, or
$n=n_i$ and $r>r_i$.  Since $a^{(1)}_{n,r} = 0$ for all $r\in\sR$
and all  $n>|V(G)|$, we may assume
(by taking a suitable subsequence) that
$n_1 = n_2 = \cdots$, and that $r_1 \le r_2\le r_3\le \cdots$.   Since
clearly $a^{(i)}_{n,r} \ge a^{(i)}_{n,r'}$ for all $n\ge 0$, all
$r\le r'$ and all $i=1,2,\ldots$, we have
$$
a^{(1)}_{n_1,r_1} > a^{(2)}_{n_1, r_1} \ge
a^{(2)}_{n_2,r_2} > a^{(3)}_{n_2, r_2}\ge
a^{(3)}_{n_3,r_3} > \cdots,$$
a contradiction. 
\end{proof}

We say that a \td\ $(T,W)$ of a graph $G$
is \dfn{minimal} if
there is no \td\ of $G$ 
of smaller size.

\begin{lemma}
\label{lem:minex}
Let $w$ be an integer, and let $G$ be a graph of tree-width at most $w$.
Then a minimal \td\ of $G$ exists, and every minimal \td\ of $G$ has width at most $w$.
\end{lemma}

\begin{proof}
The existence of a minimal \td\ follows from Lemma~\ref{wellorder}. 
If $G$ has a \td\ of width at most $w$, then every minimal \td\ has width
at most $w$, as desired.
\end{proof}

\begin{theorem}
\label{w1w6}
Let $(T,W)$ be a minimal \td\ of a graph
$G$. Then $(T,W)$ satisfies {\rm (W1)--(W6)}.
\end{theorem}

\proof That $(T,W)$ satisfies (W3) is shown in \cite{thomas90}, and that it
satisfies (W4), (W5) and (W6) is shown in \cite{OpoOxlTho}. Let us remark
that \cite{OpoOxlTho} and \cite{thomas90}
use a slightly different definition of minimality, but the proofs are
adequate, because a minimal \td\ in  our sense is minimal in the sense of
\cite{OpoOxlTho} and \cite{thomas90} as well.~\qed
\bigbreak

\begin{lemma}
\label{lem:nbrs}Let $(T,W)$ be a minimal \td\ of a graph $G$.
Then for every edge $tt'\in E(T)$ either $W_t\subseteq W_{t'}$ or
$W_{t'}\subseteq W_t$.
\end{lemma}

\begin{proof}
Assume for a contradiction that
there exists an edge $tt'\in E(T)$ such that $W_t \not\subseteq W_{t'}$ and
$W_{t'} \not \subseteq W_t$. Let $R$ be obtained from $T$ by subdividing the
edge $tt'$ and let $t''$ be the new vertex.
Let $Y_{t''} = W_t \cap W_{t'}$ and $Y_r=W_r$ for all $r\in V(T)$, 
and let $Y=(Y_r:r\in V(R))$. Then $(R,Y)$ is a \td\ of $G$ of  smaller
size than $(T,W)$,  contrary to the minimality of $(T,W)$.
\end{proof}

\begin{lemma}
\label{branchoffpath}
Let $(T,W)$ be a minimal \td\ of a 
graph $G$,
let $t\in V(T)$, let $X\subseteq W_t$,
let $B$ be a component of $T\backslash t$, let $t'$ be the neighbor of $t$ in $B$,
 let $Y=W_t\cup\bigcup_{r\in V(B)}W_{r}$,
and let $H$ be the subgraph of $G$ induced by $Y$.
 If $H\backslash X=H_1\cup H_2$, where $V(H_1)\cap V(H_2)=\emptyset$
and both of $V(H_1), V(H_2)$ intersect  $W_{t}$,
then either 
$W_{t'}-X\subseteq W_t\cap V(H_1)$
or $W_{t'}-X\subseteq W_t\cap V(H_2)$.
\end{lemma}

\begin{proof} 
%
We first prove the following claim.

\begin{claim}
\label{c:1}
Either $W_t\cap W_{t'}-X\subseteq V(H_1)$ or $W_t\cap W_{t'}-X\subseteq V(H_2)$.
\end{claim}

\noindent 
To prove the claim suppose for a contradiction that there exist vertices
$v_1\in W_t\cap W_{t'}\cap V(H_1)$ and $v_2\in W_t\cap W_{t'}\cap V(H_2)$.
Thus both $v_1$ and $v_2$ are $B$-tied, and so by (W6), which $(T,W)$
satisfies by Theorem~\ref{w1w6}, there exists a $B$-confined path $Q$ with ends $v_1$ and $v_2$.
Since $Q$ is $B$-confined, it is a subgraph of $H\backslash X$,
contrary to the fact that $V(H_1)\cap V(H_2)=\emptyset$ and
$H_1\cup H_2=H\backslash X$.
%
This proves Claim~\ref{c:1}.
\medskip

Since both of $V(H_1), V(H_2)$ intersect  $W_{t}$,
 Claim~\ref{c:1} implies that $W_t\not\subseteq W_{t'}$, and hence
$W_{t'} \subseteq W_t$ by Lemma \ref{lem:nbrs}.
By another application of Claim~\ref{c:1} we deduce that either
$W_{t'}-X\subseteq W_t\cap V(H_1)$
or $W_{t'}-X\subseteq W_t\cap V(H_2)$, as desired.
\end{proof}

%


\begin{lemma}
\label{pathsplit}Let $k\ge1$ be an  integer, let $(T,W)$ be a minimal
\td\ of a graph $G$, let $t_1,t_2\in V(T)$, let $X=W_{t_1}\cap W_{t_2}$,
let $H$ be the 
subgraph of $G$
induced by $\bigcup W_t$, the union taken over all vertices $t\in V(T)$
such that either $t\in \{t_1,t_2\}$, or for $i=1,2$ the vertex $t$ belongs
to the component of $T\backslash t_i$ containing $t_{3-i}$,
let $H\backslash X=H_1\cup H_2$, where $V(H_1)\cap V(H_2)=
\emptyset$, and assume that $|W_{t_i}\cap V(H_j)|=k$ and
$|W_t \cap V(H_i)| \ge k$  for all $i,j\in\{1,2\}$
and all $t\in T[t_1,t_2]$.
Let $t,t'$ be two adjacent vertices on the path of $T$ between $t_1$
and $t_2$.
Then  there exists an integer $i\in \{1,2\}$ such that
$W_t\cap V(H_i)=W_{t'}\cap V(H_i)$ and this set has
cardinality $k$.
\end{lemma}

\begin{proof}
%
We begin with the following claim.

\begin{claim}
\label{c:2}
For every $t\in T[t_1,t_2]$ either $|W_t \cap V(H_1)|=k$ or $|W_t \cap V(H_2)| =k$.
\end{claim}
\noindent
To prove the claim let
 $R$ be the subtree of $T$ induced by vertices $r\in V(T)$ such that either $r\in\{t_1,t_2\}$
or $r$ belongs to the component of $T\backslash \{t_1,t_2\}$ that contains neighbors of both $t_1$ and $t_2$,
let $R_1,R_2$ be two isomorphic copies of $R$, and for $r\in V(R)$  let $r_1$ and $r_2$ denote the copies
of $r$ in $R_1$ and $R_2$, respectively.
Assume for a contradiction that there is $t_0\in T[t_1,t_2]$ such that $|W_{t_0} \cap V(H_i)| > k$
for all $i \in \{1,2\}$, and choose such a vertex with $t_0\in V(R)$ and $|W_{t_0}|$ maximum.
We construct a new \td\ $(T', W')$ as follows.
The tree $T'$ is obtained from the disjoint union of $T\backslash(V(R)-\{t_1,t_2\})$, $R_1$ and $R_2$ by
identifying $t_1$ with $(t_1)_1$, $(t_2)_1$ with $(t_1)_2$ and $(t_2)_2$ with $t_2$ (here $(t_1)_2$
denotes the   copy of $t_1$ in $R_2$ and  similarly for the other three quantities).
The family $W'=(W'_t:t\in V(T'))$ is defined as follows:
$$
W'_t=\begin{cases}
W_t & \text{if }t\in V(T)-V(R)\\
(W_r\cap V(H_1))\cup (W_{t_1}\cap V(H_2)\cup X&\text{if } t=r_1 \text{ for } r\in T[t_1,t_2]\\
(W_r\cap V(H_2))\cup (W_{t_2}\cap V(H_1)\cup X&\text{if } t=r_2 \text{ for } r\in T[t_1,t_2]\\
W_r\cap V(H_1)&\text{if } t=r_1 \text{ for } r\in  V(R)-T[t_1,t_2]\\
W_r\cap V(H_2)&\text{if } t=r_2 \text{ for } r\in  V(R)- T[t_1,t_2]\\
\end{cases}
$$
Please note that the value of $W'_t$ is the same for $t=(t_2)_1$ and $t=(t_1)_2$, and  hence $W'$ is well-defined.
Since no edge of $G$ has one end in $V(H_1)$ and the other end in $V(H_2)$, it follows
that $(T',W')$ is a \td\ of $G$. 

We claim that the size of $(T',W')$ is smaller than the size of $(T,W)$. 
Indeed, let $n_0=|W_{t_0}|$, and
let $Z=\{t\in V(T'):|W'_t|\ge n_0\}$. 
Then $n_0>2k+|X|$.
We define a mapping $f:Z\to V(T)$ by $f(t)=t$ for $t\in Z-V(R_1)-V(R_2)$,
$f(r_1)=r$ for $r\in V(R)$ such that $r_1\in Z$ and $f(r_2)=r$ for $r\in V(R)$ such that $r_2\in Z$.
We remark that the vertex obtained by identifying  $(t_2)_1$ with $(t_1)_2$ does not belong to $Z$, 
and hence there is no ambiguity.
Then $Z$ and $f$ have the following properties:
\begin{itemize}
\item $|W_{f(t)}|\ge |W'_t|$ for every $t\in Z$, 
\item for $r\in V(R)$, at most one of $r_1,r_2$ belongs to $Z$, and
\item $(t_0)_1,(t_0)_2\not\in Z$
\end{itemize}
These properties follow from the assumptions  that $|W_{t_i}\cap V(H_j)|=k$ and
$|W_t \cap V(H_i)| \ge k$  for all $i,j\in\{1,2\}$
and all $t\in T[t_1,t_2]$. 
(To see  the second property assume for a contradiction that for some $r\in V(R)$
both $r_1$ and $r_2$ belong to $Z$. Then $n_0=|W_{t_0}|\ge|W_{f(r_i)}|\ge|W_{r_i}|\ge n_0$,
by the maximality of $|W_{t_0}|$ and the first property,
and  so equality holds throughout, contrary to the construction.)
It follows from the first two properties that $f$ maps injectively 
$(n,r)$-cells in $(T',W')$ to $(n,r)$-cells in $(T,W)$ for all $n\ge n_0$ and all ranks $r$.
On the other hand, the third property implies that, letting $r_1$ denote the rank of one-vertex trees,
no $(n_0,r_1)$-cell in $(T',W')$ is mapped onto the $(n_0,r_1)$-cell in $(T,W)$ with vertex-set $\{t_0\}$.
Thus  the size of $(T',W')$ is smaller than the size of $(T,W)$, contrary to the minimality of $(T,W)$.
This proves Claim~\ref{c:2}.
\medskip


Now let $t,t' \in T[t_1,t_2]$ be adjacent. By Lemma~\ref{lem:nbrs} we may assume that $W_t\subseteq W_{t'}$.
Then $W_t\cap V(H_1)\subseteq W_{t'}\cap V(H_1)$ and $W_t\cap V(H_2)\subseteq W_{t'}\cap V(H_2)$.
By Claim~\ref{c:2} we may assume that $|W_{t'}\cap V(H_1)| = k$. Given that $|W_t\cap V(H_1)| \geq k$ we have
$W_t\cap V(H_1)= W_{t'}\cap V(H_1)$ and this set has cardinality $k$, as desired.
\end{proof}


\begin{lemma}
\label{lem:triad}
Let $(T,W)$ be a minimal \td\ of a graph
$G$, let $t_1,t_2,t_3$ be a $W$-separable triad in $T$ with center
$t_0$, and let
$X,H,H_1,H_2$ and $H_3$  be as in the definition of $W$-separable triad.
Let $k=|W_{t_1}-X|/2$ and  for $i=1,2,3$ let $t_i'$ denote the neighbor
of $t_0$ in the component of $T\backslash t_0$ containing $t_i$.
Then
for all distinct $i,j\in\{1,2,3\}$, $V(H_i)\cap W_{t'_i}=\emptyset$,
$V(H_i)\cap W_{t'_j}=V(H_i)\cap W_{t_0}$, and this set has cardinality $k$.
\end{lemma}

\begin{proof}
 Let $X_3=\bigcup W_t$, the union taken over all $t\in V(T)$ that do
not belong to the component of $T\backslash t_3$ containing $t_0$.
By Lemma~\ref{branchoffpath} applied to the vertex $t_0$, the component of $T\backslash t_0$
containing $t_3$ and the subgraphs of $G$ induced by  $(V(H_1)\cup V(H_2)\cup X_3)\cap Y$
and $V(H_3)\cap Y$, where $Y$ is as in the statement of Lemma~\ref{branchoffpath} 
we deduce that $V(H_3)\cap W_{t'_3}=\emptyset$. The other two statements of the first assertion follow
by symmetry.

To prove the remaining assertions,
since $|W_{t_0}\cap V(H_1)|\ge k$ and $|W_{t_0}\cap V(H_2)|\ge k$ by the
definition of $W$-separable triad,
by Lemma~\ref{pathsplit} applied to $t_1,t_2,H_3$ and the subgraph of $G$
induced by $V(H_1)\cup V(H_2)\cup X_3$ we deduce that
$V(H_3)\cap W_{t_0}=V(H_3)\cap W_{t'_1}=V(H_3)\cap W_{t'_2}$,
and this set has cardinality $k$.
Similarly we deduce that $V(H_2)\cap W_{t_0}=V(H_2)\cap W_{t'_1}=
V(H_2)\cap W_{t'_3}$ and
$V(H_1)\cap W_{t_0}=V(H_1)\cap W_{t'_2}=V(H_1)\cap W_{t'_3}$,
and that the latter two sets also have cardinality $k$.
\end{proof}

We are finally ready to prove Theorem~\ref{w7}, which,
by Lemma~\ref{lem:minex} is implied by the following theorem.

\begin{theorem}
\label{w1w7-2}
Let $(T,W)$ be a minimal \td\ of a graph
$G$. Then $(T,W)$ satisfies {\rm (W1)--(W7)}.
\end{theorem}

\proof That $(T,W)$ satisfies (W1)--(W6) follows from Theorem~\ref{w1w6}.
Thus it remains to show that $(T,W)$ satisfies (W7).
Suppose for a contradiction that $(T,W)$ does not satisfy
(W7), and let $t_1,t_2,t_2$ be a $W$-separable triad in $T$
with center $t_0$ such that $W_{t_i}\cap W_{t_0}\subseteq X$
for every $i=1,2,3$, where $X=W_{t_1}\cap W_{t_2}\cap W_{t_3}$.
Let $H,H_1,H_2$ and $H_3$  be as in the definition of $W$-separable triad,
and for $i\in\{1,2,3\}$ let $t'_i$ denote the neighbor of $t_0$ in the component
of $T\backslash t_0$ containing $t_i$.

Let $n:=|W_{t_1}|$, let $k:=|W_{t_1}-X|/2$, let $r_1$ denote the rank of $1$-vertex trees, and let
$T_0$ denote the $(n,r_1)$-cell containing $t_0$.
By the definition of $W$-separable triad  we have $|W_{t'_i}|\ge n$ for all
$i\in\{1,2,3\}$, and hence
the degree of $t_0$ in $T_0$ is at least three and by  Lemmas~\ref{lem:triad}
and~\ref{branchoffpath} it is at most three.

Let $(T_0, P_0, T_1,P_1,\dots, T_l,P_l)$ be a spine-decomposition
of $T_0$ relative to $t_0$ with $t_0\in V(P_l)$.  
Since $P_l$ is a maximal path in $T_{l}$ we may assume that
$t'_1,t'_2\in V(P_l)$ and $t'_3\not\in V(P_l)$.

It follows from  Lemma~\ref{lem:triad} that $W_{t_3}\cap W_{t_3'}=X$.
By Lemma~\ref{pathsplit} applied to $t_3$ and $t'_3$ and $t_3'$ 
and its neighbor in $T[t_3,t_3']$ we deduce that there exists a
vertex $r_3\in T[t_3,t_3']-\{t_3'\}$ 
such that
either $V(H_1)\cap W_{t'_3}=V(H_1)\cap W_r$ for every $r\in T[r_3,t_3']$, or
$V(H_2)\cap W_{t'_3}=V(H_2)\cap W_r$ for every $r\in  T[r_3,t_3']$.
Without loss of generality we may assume the latter.
We may choose $r_3$ to be as close to $t_3$ as possible.
The fact that $W_{t_3}\cap W_{t_3'}=X$
implies that $r_3\ne t_3$.
By another application of Lemma~\ref{pathsplit}, this time to $t_3$, $t_3'$, $r_3$ and the
neighbor of $r_3$ in  $T[r_3,t_3]$, we deduce that
$|V(H_1)\cap W_{r_3}|=|V(H_2)\cap W_{r_3}|=k$.

Let $r_3'$ be the neighbor of $r_3$ in $T[r_3,t_0]$ and let
  the tree $T''$ be defined as follows: for every component $B$ of
$T\backslash T[t_0,r_3']$
not containing $t_1,t_2$ or $t_3$ let $r(B)r'(B)$ denote the edge connecting $B$ to
$T[t_0,r_3']$, where $r(B)\in V(B)$ and $r'(B)\in T[t_0,r_3']$.
By Lemma~\ref{branchoffpath} there exists an integer $i\in\{1,2,3\}$
such that $W_{r(B)}\subseteq W_{r'(B)}\cap V(H_i)$.
Let us  mention in passing that this, the choice of $r_2$ and Lemma~\ref{lem:triad}
imply that for every such component $B$, every $(n,r_1)$-cell is either a subgraph of $B$
or is disjoint from $B$.
The tree $T''$ is obtained from $T$ by, for every such component $B$ for which either $i=2$, 
or $i=3$ and $r'(B)=t_0$, deleting the edge
$r(B)r'(B)$ and adding the edge $t'_1r(B)$;
and for every such component $B$ for which  $i=1$ and $r'(B)=t_0$
deleting the edge $r(B)r'(B)$ and adding the edge $t_2'r(B)$.
Since $W_{r'(B)}\cap (V(H_2)\cup V(H_3))\subseteq  W_{t'_1}$ by the choice of $r_3$ and 
Lemma~\ref{lem:triad}; and $W_{r'(B)}\cap V(H_1)\subseteq W_{t_2'}$ by Lemma~\ref{lem:triad}
 it follows that $(T'',W)$ is a \td\ of $G$.

Let $T'$ be defined as in Lemma~\ref{spine2},
starting from the tree $T''$,
let $t'_0$ be the vertex that resulted from contracting the edges $t_0t_2'$ and $t_0t_3'$,
and let $W'=(W'_t\mid t\in V(T'))$ be defined by
$$
W'_t=\begin{cases}
W_t & \text{if }t \in V(T')-T'[r_3'',t'_0]\\
W_{r_3}\cup (V(H_3)\cap W_{t_0}) & \text{if }t=r_3'''\\
(W_{r_3}-V(H_2))\cup (V(H_3)\cap W_{t_0}) &\text{if }t=r_3''\\
W_{t_2'}&\text{if }t=t_0'\\
(W_t-V(H_2))\cup (V(H_3)\cap W_{t_0}) & \text{if } t\in T'[r_3',t'_0]-\{t'_0\}
\end{cases}
$$

We claim that $(T',W')$ is a tree decomposition  of $G$.
Indeed, since $V(H_2)\cap W_r\subseteq W_{t_0}$ for all $r\in T[r_3',t_0]$ it follows
that $(T',W')$ satisfies (W1).

To show that $(T',W')$ satisfies (W2) let $v\in V(G)$, let $Z=\{t\in V(T):v\in W_{t}\}$.
 and let $Z'=\{t\in V(T'):v\in W_{t}'\}$.
It suffices to show that $Z'$ induces a connected subset of $T'$, for this is easily seen to be
equivalent to (W2).
To that end assume first that $v\not\in W'_{t_1'}=W_{t_1'}=W_{t_0}\cap(V(H_2)\cup V(H_3)\cup X)$,
where the second equality follows from Lemma~\ref{lem:triad}.
It follows that, since $Z$ induces a subtree of $T$, that $Z'$ induces a subtree of $T'$.
We assume next that $v\in W_{t_0}\cap V(H_2)$.
The choice of $T''$ and the definition of $W'$ imply that no vertex in the component of
$T'\backslash r_3'''$ containing $t_0'$ belongs to $Z'$. Again, it follows that $Z'$ induces
a subtree of $T'$.
Finally, let $v\in W_{t_0}\cap (V(H_3)\cup X)$. Then $T[t_1',t_0']\subseteq Z'$, and it again  follows that $Z'$ induces
a subtree of $T'$. This proves our claim that $(T',W')$ is a \td.



We claim that the size of $(T',W')$ is smaller than the size of $(T,W)$. 
 Let $r$ denote the rank of $T_0$, and let $T_0'$ denote the $(n,r_1)$-cell in $(T',W')$ containing $t_0'$.
 First, by the passing remark made a few paragraphs ago, for every integer $m\ge n$ and every rank $s$,
to every $(m,s)$-cell in $(T',W')$ other than $T_0'$ there corresponds a unique
$(m,s)$-cell in $(T,W)$.
(To the $(n+1,r_1)$-cell in $(T',W')$ with vertex-set $\{r_3'''\}$ there corresponds the 
$(n+1,r_1)$-cell in $(T,W)$ with vertex-set $\{t_0\}$.)
Second, by Lemma~\ref{spine2} the rank of $T_0$ is strictly larger than the rank of $T_0'$.
Thus no $(n,r)$-cell in $(T',W')$ corresponds to $T_0$.
%
It follows that $(T',W')$ is a tree-decomposition of $G$ of
smaller size, contrary to the minimality of $(T,W)$.~\qed

%% file: image10.tikz.tex
\definecolor{qqqqff}{rgb}{0.3333333333333333,0.3333333333333333,0.3333333333333333}
\definecolor{cqcqcq}{rgb}{0.7529411764705882,0.7529411764705882,0.7529411764705882}
\begin{tikzpicture}[scale=1.4,line cap=round,line join=round,>=triangle 45,x=1.0cm,y=1.0cm]
\draw [rotate around={0.:(1.8574259450964128,2.7850457542236255)}] (1.8574259450964128,2.7850457542236255) ellipse (0.5575559570832301cm and 0.2011682014609032cm);
\draw [rotate around={0.:(0.9527563101048637,1.8780413017362978)}] (0.9527563101048637,1.8780413017362978) ellipse (0.5575559570832466cm and 0.20116820146090852cm);
\draw [rotate around={0.:(2.7639154399193266,1.8477163657537845)}] (2.7639154399193266,1.8477163657537845) ellipse (0.5575559570832466cm and 0.20116820146090852cm);
\draw [rotate around={0.:(2.0420059982729493,0.7795470836370137)}] (2.0420059982729493,0.7795470836370137) ellipse (0.557555957083252cm and 0.20116820146091097cm);
\draw [rotate around={0.:(3.652763252199438,0.7646228021339283)}] (3.652763252199438,0.7646228021339283) ellipse (0.5575559570832573cm and 0.20116820146091305cm);
\draw (1.6374259450964146,2.825045754223626)-- (0.7348160866717592,1.8867937655879032);
\draw (2.127523746118338,2.817084874632396)-- (3.055472796308959,1.8417947430774735);
\draw (2.5638717327782485,1.8295047164892058)-- (1.838760164070451,0.7848524564864476);
\draw (3.055472796308959,1.8417947430774735)-- (2.2689110946598223,0.7725624298981799);
\draw (2.4561193806564767,3.0595536307934093) node[anchor=north west] {$X_{\eta(v_1)}$};
\draw (2.640446668799136,2.662541010178453) node[anchor=north west] {$X_i(u_0)$};
\draw (0.3718031224279504,2.7200039436839933) node[anchor=north west] {$P_j(u_0)$};
\draw (2.597909602304676,1.3580709710150247) node[anchor=north west] {$R_1$};
\draw (1.7188102280858415,1.6983675029707017) node[anchor=north west] {$R_3$};
\draw [shift={(1.8194943607306575,1.2466210839566871)}] plot[domain=0.6643188953948307:2.44440039458903,variable=\t]({1.*0.9454369366016867*cos(\t r)+0.*0.9454369366016867*sin(\t r)},{0.*0.9454369366016867*cos(\t r)+1.*0.9454369366016867*sin(\t r)});
\draw [shift={(2.048131603910006,0.6384320417562718)}] plot[domain=-3.751871011414078:0.5459386178874661,variable=\t]({1.*0.2554903866493711*cos(\t r)+0.*0.2554903866493711*sin(\t r)},{0.*0.2554903866493711*cos(\t r)+1.*0.2554903866493711*sin(\t r)});
\draw (2.144180893030439,0.5073296411258326) node[anchor=north west] {$R_2$};
\draw (1.5762731615913817,2.6341829658488134) node[anchor=north west] {$Y_i(u_0)$};
\begin{scriptsize}
\draw [fill=qqqqff] (1.6374259450964146,2.825045754223626) circle (2.5pt);
\draw [fill=qqqqff] (0.7348160866717592,1.8867937655879032) circle (2.5pt);
\draw [fill=qqqqff] (2.127523746118338,2.817084874632396) circle (2.5pt);
\draw [fill=qqqqff] (3.055472796308959,1.8417947430774735) circle (2.5pt);
\draw [fill=qqqqff] (2.5638717327782485,1.8295047164892058) circle (2.5pt);
\draw [fill=qqqqff] (1.838760164070451,0.7848524564864476) circle (2.5pt);
\draw [fill=qqqqff] (2.2689110946598223,0.7725624298981799) circle (2.5pt);
\draw [fill=qqqqff] (1.0971800173576023,1.8515593465316187) circle (2.5pt);
\end{scriptsize}
\end{tikzpicture}

%% file: image11.tikz.tex
\definecolor{qqqqff}{rgb}{0.3333333333333333,0.3333333333333333,0.3333333333333333}
\definecolor{cqcqcq}{rgb}{0.7529411764705882,0.7529411764705882,0.7529411764705882}
\begin{tikzpicture}[scale=1.4,line cap=round,line join=round,>=triangle 45,x=1.0cm,y=1.0cm]
\draw [rotate around={0.0:(1.8574259450964128,2.7850457542236255)}] (1.8574259450964128,2.7850457542236255) ellipse (0.5575559570832301cm and 0.2011682014609032cm);
\draw [rotate around={0.0:(0.9527563101048637,1.8780413017362978)}] (0.9527563101048637,1.8780413017362978) ellipse (0.5575559570832466cm and 0.20116820146090852cm);
\draw [rotate around={0.0:(2.7639154399193266,1.8477163657537845)}] (2.7639154399193266,1.8477163657537845) ellipse (0.5575559570832466cm and 0.20116820146090852cm);
\draw [rotate around={0.0:(2.0420059982729493,0.7795470836370137)}] (2.0420059982729493,0.7795470836370137) ellipse (0.557555957083252cm and 0.20116820146091097cm);
\draw [rotate around={0.0:(3.652763252199438,0.7646228021339283)}] (3.652763252199438,0.7646228021339283) ellipse (0.5575559570832573cm and 0.20116820146091305cm);
\draw (1.6374259450964146,2.825045754223626)-- (0.7348160866717592,1.8867937655879032);
\draw (2.127523746118338,2.817084874632396)-- (3.055472796308959,1.8417947430774735);
\draw (3.055472796308959,1.8417947430774735)-- (2.2689110946598223,0.7725624298981799);
\draw (2.452483070846853,3.0110813043153994) node[anchor=north west] {$X_{\eta(v_1)}$};
\draw (2.640537711349196,2.605279185336659) node[anchor=north west] {$X_i(u_0)$};
\draw (0.3739844126630567,2.675586347362605) node[anchor=north west] {$P_j(u_0)$};
\draw (2.1722782961628886,1.6155179195348528) node[anchor=north west] {$R_1$};
\draw [shift={(1.8194943607306575,1.2466210839566871)}] plot[domain=0.6643188953948307:2.44440039458903,variable=\t]({1.0*0.9454369366016867*cos(\t r)+-0.0*0.9454369366016867*sin(\t r)},{0.0*0.9454369366016867*cos(\t r)+1.0*0.9454369366016867*sin(\t r)});
\draw [shift={(2.048131603910006,0.6384320417562718)}] plot[domain=-3.751871011414078:0.5459386178874661,variable=\t]({1.0*0.2554903866493711*cos(\t r)+-0.0*0.2554903866493711*sin(\t r)},{0.0*0.2554903866493711*cos(\t r)+1.0*0.2554903866493711*sin(\t r)});
\draw (2.0555546911063104,0.4377020132307036) node[anchor=north west] {$R_2$};
\draw (3.0555546911063104,1.4377020132307036) node[anchor=north west] {$R_3$};
\draw (1.5804692835214428,2.5854839600206225) node[anchor=north west] {$Y_i(u_0)$};
\draw [shift={(2.6615636253225765,0.083509338354617)}] plot[domain=0.6848212459425996:2.43571715830806,variable=\t]({1.0*1.0749097337369704*cos(\t r)+-0.0*1.0749097337369704*sin(\t r)},{0.0*1.0749097337369704*cos(\t r)+1.0*1.0749097337369704*sin(\t r)});
\draw (2.5638717327782485,1.8295047164892058)-- (3.494115748357463,0.7634269054997732);
\begin{scriptsize}
\draw [fill=qqqqff] (4.172763252199429,0.7646228021339283) circle (2.5pt);
\draw [fill=qqqqff] (1.6374259450964146,2.825045754223626) circle (2.5pt);
\draw [fill=qqqqff] (0.7348160866717592,1.8867937655879032) circle (2.5pt);
\draw [fill=qqqqff] (2.127523746118338,2.817084874632396) circle (2.5pt);
\draw [fill=qqqqff] (3.055472796308959,1.8417947430774735) circle (2.5pt);
\draw [fill=qqqqff] (2.5638717327782485,1.8295047164892058) circle (2.5pt);
\draw [fill=qqqqff] (1.838760164070451,0.7848524564864476) circle (2.5pt);
\draw [fill=qqqqff] (2.2689110946598223,0.7725624298981799) circle (2.5pt);
\draw [fill=qqqqff] (1.0971800173576023,1.8515593465316187) circle (2.5pt);
\draw [fill=qqqqff] (3.494115748357463,0.7634269054997732) circle (2.5pt);
\end{scriptsize}
\end{tikzpicture}

%% file: image13.tikz.tex
\definecolor{qqqqff}{rgb}{0.3333333333333333,0.3333333333333333,0.3333333333333333}
\begin{tikzpicture}[scale=1.4,line cap=round,line join=round,>=triangle 45,x=1.0cm,y=1.0cm]
\clip(-2.05598612148122,-1.7511229607967027) rectangle (10.94054048473452,3.006101850147045);
\draw [rotate around={0.:(1.8574259450964128,2.7850457542236255)}] (1.8574259450964128,2.7850457542236255) ellipse (0.5575559570832301cm and 0.2011682014609032cm);
\draw [rotate around={0.:(0.9527563101048637,1.8780413017362978)}] (0.9527563101048637,1.8780413017362978) ellipse (0.5575559570832466cm and 0.20116820146090852cm);
\draw [rotate around={0.:(2.7639154399193266,1.8477163657537845)}] (2.7639154399193266,1.8477163657537845) ellipse (0.5575559570832466cm and 0.20116820146090852cm);
\draw [rotate around={0.:(3.4941743671537235,0.9381472310840802)}] (3.4941743671537235,0.9381472310840802) ellipse (0.5575559570832597cm and 0.20116820146091385cm);
\draw (1.6374259450964146,2.825045754223626)-- (0.7348160866717592,1.8867937655879032);
\draw (2.127523746118338,2.817084874632396)-- (3.1249622632109437,1.8614589794739969);
\draw (2.4554368145906063,3.089641380856199) node[anchor=north west] {$X_{\eta(v_1)}$};
\draw (2.5431010347833626,2.6868153491861926) node[anchor=north west] {$X_i(u_0)$};
\draw (0.37772580531366046,2.6384762253857916) node[anchor=north west] {$P_j(u_0)$};
\draw [shift={(1.8194943607306575,1.2466210839566871)}] plot[domain=0.7057288752930431:2.417996438474012,variable=\t]({1.*0.8797138490345054*cos(\t r)+0.*0.8797138490345054*sin(\t r)},{0.*0.8797138490345054*cos(\t r)+1.*0.8797138490345054*sin(\t r)});
\draw (1.6779707595063298,2.504589266652592) node[anchor=north west] {$Y_i(u_0)$};
\draw (2.489078892078955,1.8171936135812614)-- (3.1980078738131668,0.9165135049280175);
\draw [rotate around={0.:(2.8922318476463027,-0.14952232506123514)}] (2.8922318476463027,-0.14952232506123514) ellipse (0.5575559570832355cm and 0.20116820146090453cm);
\draw [rotate around={0.:(4.559839449215463,-0.17948716665584855)}] (4.559839449215463,-0.17948716665584855) ellipse (0.5575559570832155cm and 0.20116820146089792cm);
\draw [rotate around={0.:(3.883213957274235,-1.1259011294376284)}] (3.883213957274235,-1.1259011294376284) ellipse (0.5575559570832701cm and 0.20116820146091727cm);
\draw [rotate around={0.:(5.377299276574322,-1.118146734395608)}] (5.377299276574322,-1.118146734395608) ellipse (0.5575559570832647cm and 0.20116820146091569cm);
\draw (3.1249622632109437,1.8614589794739969)-- (3.8136837671123778,0.9774609268976735);
\draw (3.1980078738131668,0.9165135049280175)-- (2.7996293546195594,-0.110885893989041);
\draw (3.8136837671123778,0.9774609268976735)-- (4.733575984159435,-0.1631265413916033);
\draw (4.733575984159435,-0.1631265413916033)-- (5.4796303019220085,-1.129578518339006);
\draw (5.117468011746001,-1.155698842040287)-- (4.364170448179908,-0.19795363965997817);
\draw (4.364170448179908,-0.19795363965997817)-- (3.5022241975610124,0.9339270540622049);
\draw (3.5022241975610124,0.9339270540622049)-- (3.1183321699744453,-0.13700621769032215);
\draw [shift={(2.937251024886442,-0.27631461076382163)}] plot[domain=-4.018491769346692:0.6557483755750777,variable=\t]({1.*0.21518918295489942*cos(\t r)+0.*0.21518918295489942*sin(\t r)},{0.*0.21518918295489942*cos(\t r)+1.*0.21518918295489942*sin(\t r)});
\draw [shift={(5.3057924026375245,-1.2950072351137865)}] plot[domain=-3.7784857215343823:0.76061688328206,variable=\t]({1.*0.23424966293547556*cos(\t r)+0.*0.23424966293547556*sin(\t r)},{0.*0.23424966293547556*cos(\t r)+1.*0.23424966293547556*sin(\t r)});
\draw (2.3699102090966136,1.6072415843105752) node[anchor=north west] {$R_1$};
\draw (3.4593001174097774,0.5115547781681578) node[anchor=north west] {$R_2$};
\draw (4.08462855779529,0.833815603504163) node[anchor=north west] {$R_3$};
\draw (2.7905514935062428,-0.42632336770545556) node[anchor=north west] {$R_4$};
\draw (5.3765680073855735,-1.3414449675138718) node[anchor=north west] {$R_5$};
\draw (3.5814221357471445,-1.2125406373794696) node[anchor=north west] {$X_{\eta(v_8)}$};
\draw (0.4581533282534132,1.7039198319113769) node[anchor=north west] {$X_{\eta(v_2)}$};
\begin{scriptsize}
\draw [fill=qqqqff] (1.6374259450964146,2.825045754223626) circle (2.5pt);
\draw [fill=qqqqff] (0.7348160866717592,1.8867937655879032) circle (2.5pt);
\draw [fill=qqqqff] (2.127523746118338,2.817084874632396) circle (2.5pt);

\draw [fill=qqqqff] (3.1249622632109437,1.8614589794739969) circle (2.5pt);
\draw [fill=qqqqff] (2.489078892078955,1.8171936135812614) circle (2.5pt);
\draw [fill=qqqqff] (1.1495618246053554,1.8384734225828168) circle (2.5pt);
\draw [fill=qqqqff] (3.1980078738131668,0.9165135049280175) circle (2.5pt);

\draw [fill=qqqqff] (3.8136837671123778,0.9774609268976735) circle (2.5pt);
\draw [fill=qqqqff] (2.7996293546195594,-0.110885893989041) circle (2.5pt);
\draw [fill=qqqqff] (4.733575984159435,-0.1631265413916033) circle (2.5pt);
\draw [fill=qqqqff] (5.4796303019220085,-1.129578518339006) circle (2.5pt);
\draw [fill=qqqqff] (5.117468011746001,-1.155698842040287) circle (2.5pt);
\draw [fill=qqqqff] (4.364170448179908,-0.19795363965997817) circle (2.5pt);
\draw [fill=qqqqff] (3.5022241975610124,0.9339270540622049) circle (2.5pt);
\draw [fill=qqqqff] (3.1183321699744453,-0.13700621769032215) circle (2.5pt);
\end{scriptsize}
\end{tikzpicture}

%% file: image14.tikz.tex
\definecolor{qqqqff}{rgb}{0.3333333333333333,0.3333333333333333,0.3333333333333333}
\begin{tikzpicture}[scale=1.4,line cap=round,line join=round,>=triangle 45,x=1.0cm,y=1.0cm]
\clip(-2.1989750023015024,-1.06497277009617) rectangle (8.95405802218279,3.137615617282007);
\draw [rotate around={0.:(1.8574259450964128,2.7850457542236255)}] (1.8574259450964128,2.7850457542236255) ellipse (0.5575559570832301cm and 0.2011682014609032cm);
\draw [rotate around={0.:(0.9527563101048637,1.8780413017362978)}] (0.9527563101048637,1.8780413017362978) ellipse (0.5575559570832466cm and 0.20116820146090852cm);
\draw [rotate around={0.:(2.7639154399193266,1.8477163657537845)}] (2.7639154399193266,1.8477163657537845) ellipse (0.5575559570832466cm and 0.20116820146090852cm);
\draw [rotate around={0.:(3.5358810695354466,0.929791611116991)}] (3.5358810695354466,0.929791611116991) ellipse (0.5575559570832398cm and 0.20116820146090658cm);
\draw (2.127523746118338,2.817084874632396)-- (3.0646947029112397,1.8752196219248107);
\draw (2.4533143784889657,3.04022664868325) node[anchor=north west] {$X_{\eta(v_1)}$};
\draw (2.64023671968144,2.640776339504564) node[anchor=north west] {$X_j(u_0)$};
\draw (1.2848558646976832,2.615810695180896) node[anchor=north west] {$Y_j(u_0)$};
\draw (2.598969859648616,1.8417971420564534)-- (3.3357882683924682,0.9477458055778964);
\draw [rotate around={0.:(2.8852807305826853,-0.01583240558780602)}] (2.8852807305826853,-0.01583240558780602) ellipse (0.5575559570832382cm and 0.2011682014609058cm);
\draw [rotate around={0.:(4.427768225006651,0.004336472620116452)}] (4.427768225006651,0.004336472620116452) ellipse (0.5575559570832196cm and 0.2011682014608997cm);
\draw (3.0646947029112397,1.8752196219248107)-- (3.815415345782334,0.9811682854462537);
\draw (3.815415345782334,0.9811682854462537)-- (4.600891573971535,0.011916369263892774);
\draw (0.449091497925215,1.6546333887196836) node[anchor=north west] {$X_{\eta(v_2)}$};
\draw (3.3357882683924682,0.9477458055778964)-- (4.121264496581669,-0.004794870670285864);
\draw (2.835307839811738,1.8501527620235427)-- (1.889955919159249,2.8026934382717252);
\draw [rotate around={0.:(2.2500802266959585,0.9226789456766284)}] (2.2500802266959585,0.9226789456766284) ellipse (0.5575559570832523cm and 0.20116820146091127cm);
\draw (2.835307839811738,1.8501527620235427)-- (2.459947518376191,0.9059677057424498);
\draw (2.0220271433680526,0.9310345656437178)-- (2.390436347739979,1.8417971420564534);
\draw [shift={(1.78568916320493,1.14828068478804)}] plot[domain=0.8536672874933118:2.218305224992563,variable=\t]({1.*0.9201544618732224*cos(\t r)+0.*0.9201544618732224*sin(\t r)},{0.*0.9201544618732224*cos(\t r)+1.*0.9201544618732224*sin(\t r)});
\draw [shift={(1.6744712901869903,0.8224115060715567)}] plot[domain=0.834172576332828:2.225665532671888,variable=\t]({1.*1.3761702946179752*cos(\t r)+0.*1.3761702946179752*sin(\t r)},{0.*1.3761702946179752*cos(\t r)+1.*1.3761702946179752*sin(\t r)});
\draw (2.5779292726172818,-0.15537582474623618) node[anchor=north west] {$X_{\eta(v_6)}$};
\draw [shift={(2.2144691647430266,0.7433873709058068)}] plot[domain=-3.914376483006601:0.5849734910252258,variable=\t]({1.*0.2687850466151223*cos(\t r)+0.*0.2687850466151223*sin(\t r)},{0.*0.2687850466151223*cos(\t r)+1.*0.2687850466151223*sin(\t r)});
\draw [shift={(4.384845237477866,-0.16785864690807012)}] plot[domain=-3.695611407338282:0.6940177865128512,variable=\t]({1.*0.30994290134202834*cos(\t r)+0.*0.30994290134202834*sin(\t r)},{0.*0.30994290134202834*cos(\t r)+1.*0.30994290134202834*sin(\t r)});
\draw (4.077153747264006,0.8181468811714817) node[anchor=north west] {$Q_{l_j}(u_2)$};
\begin{scriptsize}
\draw [fill=qqqqff] (1.889955919159249,2.8026934382717252) circle (2.5pt);
\draw [fill=qqqqff] (2.127523746118338,2.817084874632396) circle (2.5pt);
\draw [fill=qqqqff] (3.0646947029112397,1.8752196219248107) circle (2.5pt);

\draw [fill=qqqqff] (2.598969859648616,1.8417971420564534) circle (2.5pt);
\draw [fill=qqqqff] (3.3357882683924682,0.9477458055778964) circle (2.5pt);

\draw [fill=qqqqff] (3.815415345782334,0.9811682854462537) circle (2.5pt);
\draw [fill=qqqqff] (4.600891573971535,0.011916369263892774) circle (2.5pt);
\draw [fill=qqqqff] (4.121264496581669,-0.004794870670285864) circle (2.5pt);

\draw [fill=qqqqff] (2.835307839811738,1.8501527620235427) circle (2.5pt);

\draw [fill=qqqqff] (2.459947518376191,0.9059677057424498) circle (2.5pt);
\draw [fill=qqqqff] (2.0220271433680526,0.9310345656437178) circle (2.5pt);
\draw [fill=qqqqff] (2.390436347739979,1.8417971420564534) circle (2.5pt);

\draw [fill=qqqqff] (1.2295997981152305,1.8835752418919) circle (2.5pt);
\draw [fill=qqqqff] (0.8403372425524408,1.908642101793168) circle (2.5pt);
\end{scriptsize}
\end{tikzpicture}

%% file: image12.tikz.tex
\definecolor{xdxdff}{rgb}{0.6588235294117647,0.6588235294117647,0.6588235294117647}
\definecolor{qqqqff}{rgb}{0.3333333333333333,0.3333333333333333,0.3333333333333333}
\definecolor{cqcqcq}{rgb}{0.7529411764705882,0.7529411764705882,0.7529411764705882}
\begin{tikzpicture}[scale=1.4,line cap=round,line join=round,>=triangle 45,x=1.0cm,y=1.0cm]
\draw [rotate around={0.:(2.6393373450798454,3.141359809912277)}] (2.6393373450798454,3.141359809912277) ellipse (0.5575559570832693cm and 0.20116820146091732cm);
\draw [rotate around={0.:(1.3684560417416203,2.0166078789485504)}] (1.3684560417416203,2.0166078789485504) ellipse (0.5575559570832501cm and 0.2011682014609105cm);
\draw [rotate around={0.:(3.866834343572685,1.9883086959413339)}] (3.866834343572685,1.9883086959413339) ellipse (0.5575559570832448cm and 0.20116820146090844cm);
\draw [rotate around={0.:(3.2693099678671937,0.838932759585122)}] (3.2693099678671937,0.838932759585122) ellipse (0.5575559570832497cm and 0.20116820146091022cm);
\draw [rotate around={0.:(4.57426630884901,0.8395904506320724)}] (4.57426630884901,0.8395904506320724) ellipse (0.55755595708321cm and 0.20116820146089584cm);
\draw (3.2145992455142456,3.4069858106361237) node[anchor=north west] {$X_{\eta(v_1)}$};
\draw [rotate around={0.:(0.6366464835545873,0.8336065195514261)}] (0.6366464835545873,0.8336065195514261) ellipse (0.5575559570832497cm and 0.20116820146091022cm);
\draw [rotate around={0.:(1.962167215146671,0.8354724542891597)}] (1.962167215146671,0.8354724542891597) ellipse (0.5575559570832578cm and 0.20116820146091355cm);
\draw (2.4524830708468537,3.149647881527654)-- (1.1558958126464853,2.0213200385135943);
\draw (1.1558958126464853,2.0213200385135943)-- (0.38388202532107507,0.8534017448674627);
\draw (2.9374660910897394,3.149647881527654)-- (4.154872448025963,2.0213200385135943);
\draw (4.154872448025963,2.0213200385135943)-- (4.857602946745247,0.8435041322094445);
\draw (1.6705716708634253,2.0213200385135943)-- (2.234735592370456,0.8336065195514265);
\draw (3.6,2.)-- (3.0562374429859562,0.8336065195514265);
\draw [shift={(1.3142576151747745,0.7049375549971917)}] plot[domain=0.3016190462662392:2.8387077852148224,variable=\t]({1.*0.4664496086233137*cos(\t r)+0.*0.4664496086233137*sin(\t r)},{0.*0.4664496086233137*cos(\t r)+1.*0.4664496086233137*sin(\t r)});
\draw [shift={(3.9470225822075835,0.6554494917071013)}] plot[domain=0.4076315054576267:2.73670086730471,variable=\t]({1.*0.4743631576606497*cos(\t r)+0.*0.4743631576606497*sin(\t r)},{0.*0.4743631576606497*cos(\t r)+1.*0.4743631576606497*sin(\t r)});
\draw [shift={(2.630640098691179,1.1008420613179142)}] plot[domain=0.7478448913384683:2.377244012329392,variable=\t]({1.*1.3221738225212543*cos(\t r)+0.*1.3221738225212543*sin(\t r)},{0.*1.3221738225212543*cos(\t r)+1.*1.3221738225212543*sin(\t r)});
\draw (0.8168085893850047,1.5083419315282782)-- (1.2359876749503687,1.164773453815575);
\draw (4.47754585742123,1.4805012255834995)-- (4.101375687907944,1.1039977322637173);
\draw [shift={(1.997192888578022,0.7346303929712459)}] plot[domain=-3.5713549326864804:0.3947911196997607,variable=\t]({1.*0.26130447225550374*cos(\t r)+0.*0.26130447225550374*sin(\t r)},{0.*0.26130447225550374*cos(\t r)+1.*0.26130447225550374*sin(\t r)});
\draw [shift={(3.3220060242952813,0.7346303929712459)}] plot[domain=-3.57135493268648:0.4297622790966883,variable=\t]({1.*0.2613044722555041*cos(\t r)+0.*0.2613044722555041*sin(\t r)},{0.*0.2613044722555041*cos(\t r)+1.*0.2613044722555041*sin(\t r)});
\draw (0.431322341771527,1.7045964334570163) node[anchor=north west] {$x$};
\draw (1.2746671645427021,1.4472585043485466) node[anchor=north west] {$y$};
\draw (4.590367404978759,1.7243916587730523) node[anchor=north west] {$z$};
\draw (3.863746187264162,1.4592038638462033) node[anchor=north west] {$w$};
\draw (2.1060666278162206,0.6168829144948484) node[anchor=north west] {$P_1$};
\draw (3.039855468268849,0.53770201323070385) node[anchor=north west] {$P_2$};
\draw (3.55111807588686,2.902207565077202) node[anchor=north west] {$R_j(u_0)$};
\draw (1.0043600025167563,3.020978916973419) node[anchor=north west] {$R_i(u_0)$};
\draw (2.06442843034451,2.852719501787112) node[anchor=north west] {$R_{ij}(u_0)$};
\begin{scriptsize}
\draw [fill=qqqqff] (2.4524830708468537,3.149647881527654) circle (2.5pt);
\draw [fill=qqqqff] (1.1558958126464853,2.0213200385135943) circle (2.5pt);
\draw [fill=qqqqff] (0.38388202532107507,0.8534017448674627) circle (2.5pt);
\draw [fill=qqqqff] (2.9374660910897394,3.149647881527654) circle (2.5pt);
\draw [fill=qqqqff] (4.154872448025963,2.0213200385135943) circle (2.5pt);
\draw [fill=qqqqff] (4.857602946745247,0.8435041322094445) circle (2.5pt);
\draw [fill=qqqqff] (1.6705716708634253,2.0213200385135943) circle (2.5pt);
\draw [fill=qqqqff] (2.234735592370456,0.8336065195514265) circle (2.5pt);
\draw [fill=qqqqff] (3.6,2.) circle (2.5pt);
\draw [fill=qqqqff] (3.0562374429859562,0.8336065195514265) circle (2.5pt);
\draw [fill=qqqqff] (1.7596501847855879,0.8435041322094445) circle (2.5pt);
\draw [fill=qqqqff] (0.8391722075899067,0.8534017448674627) circle (2.5pt);
\draw [fill=qqqqff] (4.382517539160379,0.8435041322094445) circle (2.5pt);
\draw [fill=qqqqff] (3.5313228505708243,0.8336065195514265) circle (2.5pt);
\draw [fill=xdxdff] (0.8168085893850047,1.5083419315282782) circle (2.5pt);
\draw [fill=xdxdff] (1.2359876749503687,1.164773453815575) circle (2.5pt);
\draw [fill=xdxdff] (4.47754585742123,1.4805012255834995) circle (2.5pt);
\draw [fill=xdxdff] (4.101375687907944,1.1039977322637173) circle (2.5pt);
\end{scriptsize}
\end{tikzpicture}